\documentclass[ejs]{imsart}

\RequirePackage[OT1]{fontenc}
\RequirePackage{amsthm,amsmath}
\RequirePackage[numbers]{natbib}
\RequirePackage[colorlinks,citecolor=blue,urlcolor=blue]{hyperref}

\usepackage{a4}
\usepackage{graphicx}
\usepackage{parskip}

\usepackage{amsfonts}
\usepackage{natbib}
\usepackage{algorithm}
\usepackage{algorithmic}

\usepackage{latexsym}
\usepackage{latexsym,amsmath,amssymb,amsthm,epic,eepic,multirow}

\newtheorem{theorem}{Theorem}[section]

\newtheorem{lemma}{Lemma}[section]
\newtheorem{corollary}{Corollary}[section]
\newtheorem{remark}{Remark}[section]

 \newcommand{\vvert} {| \hskip-0.15em |  \hskip-0.15em |}
\newcommand{\R}{\mathbb{R}}
\newcommand{\Nat}{\mathbb{N}}

\newcommand{\PP} {{  \rm I\hskip-0.22em P}}

\newcommand{\EE} {{\rm I\hskip-0.48em E}}

\usepackage{natbib}

\pubyear{2006}
\volume{0}
\issue{0}
\firstpage{1}
\lastpage{8}

\startlocaldefs
\numberwithin{equation}{section}
\theoremstyle{plain}

\endlocaldefs

\begin{document}

\begin{frontmatter}
\title{Worst possible sub-directions in high-dimensional models}
\runtitle{Worst possible sub-directions}
%\thankstext{T1}{Footnote to the title with the `thankstext' command.}

\begin{aug}
\author{\fnms{Sara} \snm{van de Geer}\ead[label=e1]{geer@stat.math.ethz.ch}}
%\and
%\author{\fnms{Second} \snm{Author}\thanksref{t3}\ead[label=e2]{second@somewhere.com}}

\address{Seminar for Statistics\\
ETH Z\"urich\\
R\"amistrasse 101\\
8092 Z\"urich\\
Switzerland\\ 
\printead{e1}}

%\author{\fnms{Third} \snm{Author}
%\ead[label=e3]{third@somewhere.com}
%\ead[label=u1,url]{www.foo.com}}

%\address{Address of the Third author\\
%usually few lines long\\
%usually few lines long\\
%\printead{e3}\\
%\printead{u1}}
%
%\thankstext{t1}{Some comment}
%\thankstext{t2}{First supporter of the project}
%\thankstext{t3}{Second supporter of the project}
%\runauthor{F. Author et al.}

\affiliation{ETH Z\"urich}

\end{aug}

\begin{abstract}
We examine the rate of convergence of the Lasso estimator of lower dimensional
components of the high-dimensional parameter. Under bounds on the $\ell_1$-norm
on the worst possible sub-direction these rates are of order $\sqrt {|J| \log p / n }$
where $p$ is the total number of parameters, $J \subset \{ 1 , \ldots , p \}$ represents
a subset of the parameters and
$n$ is the number of observations. We also derive rates in sup-norm
in terms of the rate of convergence in $\ell_1$-norm. The irrepresentable condition
on a set $J$ requires that the $\ell_1$-norm of the worst possible sub-direction is 
sufficiently smaller than one.
In that case sharp oracle results can be obtained. Moreover, if the
coefficients in $J$ are small  enough the Lasso will
put these coefficients to zero. This extends known results which say that
the irrepresentable condition on the inactive set
(the set where coefficients are exactly zero) implies no false positives.
We further show that by de-sparsifying one obtains fast rates in supremum
norm without conditions on the worst possible sub-direction.
The main assumption here is that approximate
sparsity is of order $o (\sqrt n / \log p )$.
The results are extended to M-estimation with $\ell_1$-penalty for generalized linear models
and exponential families for example. For the graphical Lasso this leads to an
extension of known results to the case where the precision matrix is only approximately
sparse.  The bounds we provide are non-asymptotic but we also present asymptotic
formulations for ease of interpretation. 
\end{abstract}



%\tableofcontents
\end{frontmatter}

\thispagestyle{plain}

\section{Introduction}
We consider estimation bounds for parameters of interest in high-dimensional models.
We apply the M-estimation procedure with
$\ell_1$-penalty and show that under certain conditions one-dimensional parameters
can be estimated with rate $\sqrt {\log p /n}$, where $p$ is
the total number of parameters and $n$ is the number of
observations.  More generally, for a subset $J \subset \{ 1 , \ldots , p \}$
the group of parameters with index in $J$ can be estimated with rate
$\sqrt { |J | \log p / n }$ in $\ell_2$-norm. For this to happen, 
the ``worst possible sub-direction" is required to have a bounded
$\ell_1$-norm. If this $\ell_1$-norm is less than one,
we obtain oracle rates $\sqrt {| J \cap S| \log p / n}$
where $S$ is the set of active parameters or a sparse approximation thereof.
Taking $J$ to be the set $S^c$ gives variable selection results under an irrepresentable condition.

By de-sparsifying one obtains fast rates (of order $1/\sqrt n$ under certain conditions) for
one-dimensional parameters without conditions
on the $\ell_1$-norm of the worst possible sub-direction.
The de-sparsified estimator can moreover be used for
the construction of asymptotic confidence intervals
for parameters of interest. 
This study is an intermediate step towards this end.
We investigate the rates and conditions for 
remainder terms to be negligible.
Global convergence (e.g.\ in $\ell_1$-norm) is generally sufficient for the latter.
However, in high-dimensional
models which are not very sparse, global convergence 
does not happen. For example, when estimating a
$p \times p$ precision matrix (where there
are actually $p(p-1)/2$ parameters), the global
rate in $\ell_1$-norm will not be faster than
$p \sqrt{ \log p / n }$. To handle such cases,
we show that the irrepresentable condition on the set of small coefficients yield
rates in the sup-norm.
 
 \subsection{Related work}
 
 The motivation of this study is founded in \cite{Jankova2014} where asymptotic
 confidence intervals for the elements of a precision matrix are studied, based
 on the graphical Lasso. This work uses results from \cite{ravikumar2011high},
 which in turn relies on irrepresentable conditions implying that with high probability there are no false positives.
 In this paper we extend such a result to the case
 where the model is only approximately sparse which is
 possibly more appropriate in the context
 of confidence intervals and testing.
 
 The literature on a semi-parametric approach to confidence intervals and testing 
 in high dimensions is expanding quickly.
 An important reference is \cite{zhangzhang11} and further work can be found
 in \cite{jamo13}, \cite{jamo13b}, \cite{vdgetal13} and the papers
 \cite{Belloni2013b}, 
 \cite{Belloni2013a},  and \cite{beletal11}. Our work
 presents rates in sup-norm for Lasso estimators and is in that
 aspect related to 
 \cite{lounici2008sup} although our conditions are based
 on worst possible sub-directions instead of incoherence.
 Also related is \cite{wainwright2009sharp} but our work
 does not rely on irrepresentable conditions, i.e.,  the
 $\ell_1$-norm of the worst possible sub-direction is allowed
 to be larger than one (but if it is smaller we reproduce variable
 selection results). 
Irrepresentable conditions for variable selection were introduced
in \cite{mebu06} and \cite{ZY07}. Our formulation shows these
are conditions on 
worst possible sub-directions. We moreover extend the situation to
models which are only approximately sparse.

 \subsection{Organization of the paper}
 The paper is organized as follows. In Section \ref{fixed-design.section}
 we consider the linear model with fixed design and the Lasso estimator.
 We derive in Subsection \ref{single.section} rates for a single coefficient and in 
 Subsection \ref{infinity.section} rates in (weighted) sup-norm.
 We consider de-sparsifying the Lasso in Subsection \ref{de-sparsity.section},
 leading to improved rates. We also discuss 
 thresholding yielding a re-sparsified estimator.
 The results are based on approximate worst possible sub-directions
 using a Lasso but one can also apply a Dantzig selector. This is
 discussed in Subsection 
\ref{Dantzig.section}. Subsection \ref{group.section} gives rates for groups of
variables. Sharp oracle inequalities as
well as variable selection results are derived. This leads to a further
refinement in Subsection \ref{selection.section} where we prove
that under certain irrepresentable conditions the Lasso will estimate small coefficients as being zero.
In the final part of this section,  Subsection \ref{de-group.section}, we present results for a
de-sparsified estimator of a group of variables. 

In the remainder of the paper worst possible sub-directions are
in terms of theoretical (unknown) quantities, which means they
do not immediately lead to a de-sparsifying procedure. 
We remark that de-sparsifying is nevertheless possible
(see also \cite{vdgetal13} and \cite{Jankova2014}) but a full discussion goes beyond the scope of this paper.
Section \ref{random.section} gives a result analogous to
the one of Subsection \ref{group.section} for groups of variables
for the case of random design. Here, worst possible sub-directions
are taken in terms of the population inner-product matrix.
Section \ref{general.section} studies general loss functions.
In Section \ref{rem.section} we discuss the remainder term,
for the linear model with random design (Subsection \ref{rem-linear.section}),
the generalized linear model (Subsection \ref{rem-GLM.section}) and
exponential families (Subsection \ref{rem-exponential.section}).
Then we move to Brouwer's fixed point theorem for deriving rates
for estimators defined as solution of a system of equations.
This theorem provides a way to handle the  situation where the
global rate is not fast enough to deal with the remainder term.
We apply this in Section \ref{irrepresentable.section}
to derive rates in sup-norm from the KKT conditions. Finally, we examine
 in Section \ref{precision.section} the remainder term of the graphical Lasso as an example.
 The approach there is as in 
 \cite{ravikumar2011high} but with the extension to models
 which are only approximately sparse. Section \ref{proofs.section} contains all proofs.
 
 The results in this paper are presented in a non-asymptotic form. To simplify their
 interpretation, we present a separate asymptotic formulation at various
 stages, where we assume ``standard" asymptotic scenarios. 
   
 \section{The linear model with fixed design}\label{fixed-design.section}
 
 Let $Y$ be an $n$-vector of response variables and
 $X$ a fixed $n \times p$ design matrix and consider the model
 \begin{equation}\label{linearmodel.equation}
 Y = X \beta^0 + \epsilon , 
 \end{equation} 
 where $\epsilon$ is unobservable noise and $\beta^0$ is a $p$-vector
 of unknown coefficients. 
 For a vector $v \in \R^n$ we write (with some abuse of notation) $\| v \|_n^2 := v^T v / n$.
The Lasso estimator (\cite{tibs96}) is
\begin{equation}\label{Lasso.equation}
 \hat \beta := \hat \beta (\lambda) :=  \min_{\beta \in \R^p} \biggl \{ \| Y - X \beta \|_n^2 + 2 \lambda \| \beta \|_1 \biggr \} . 
 \end{equation}
Here, $\lambda$ is a tuning parameter which may be chosen data-dependent
(e.g.\ when using the square root Lasso introduced in \cite{belloni2011square}). Typically, $\lambda$ is chosen
of order $\sqrt {\log p / n}$ and proportional to some estimate
of the noise level $\sigma_{\epsilon}:=  (\EE \| \epsilon \|_n^2)^{1/2}  $.

The estimator $\hat \beta$ satisfies the Karush-Kuhn-Tucker or KKT conditions
\begin{equation} \label{KKT.equation}
-X^T (Y- X \hat \beta )/n + \lambda \hat z =0 
\end{equation}
where $\hat z_j = {\rm sign} ( \hat \beta_j) $ if $\hat \beta_j \not=0 $ and $\| \hat z \|_{\infty} \le 1 $.
Thus, $\hat \beta^T \hat z = \| \hat \beta \|_1$ and
\begin{equation}\label{useful.equation}
Y^T (Y- X \hat \beta)/n = \| Y- X \hat \beta \|_n^2 + \lambda \| \hat \beta \|_1 . 
\end{equation} 
These equalities will play a key role in our proofs.

\subsection{Bounds for a single parameter}\label{single.section}
Let $j \in \{ 1 , \ldots , p \} $ be some index. We define (approximate) worst
possible sub-directions with help of the Lasso, where we regress $X_j$ on
the set of all other variables $X_{-j} := \{ X_k : \ k \not= j \} $ with $\ell_1$-penalty on the coefficients:
$$\hat \gamma_j := \hat \gamma_j ( \lambda_j) := \arg \min_{\gamma_j \in \R^{p-1} } \biggl \{ 
\| X_j - X_{-j} \gamma_j  \|_n^2 + 2 \lambda_j \| \gamma_j \|_1 \biggr \} . $$
We leave the choice of the tuning parameter
$\lambda_j$
free at this stage, but will indicate in Corollary \ref{oneparameter.corollary} that the square root Lasso 
gives well- scaled bounds.

Define $\hat \tau_j^2 := \| X_j - X_{-j}\hat \gamma_j  \|_n^2 $
and $\tilde \tau_j^2  := \hat \tau_j^2  + \lambda_j \| \hat \gamma_j \|_1 $.
Let $\hat c_{k,j}  := - \gamma_{k,j} $, $k \not=j $ and $\hat c_{j,j} := 1 $.
Note that $\hat \tau_j = \| X \hat c_j \|_n $. 
Inspired by semi-parametric theory (see e.g. \cite{bicketal98}), we call $\hat c_j$ the (approximate) worst possible
sub-direction for estimating $\beta_j^0$. We further write
$\hat \theta_j := \hat c_j / \tilde \tau_j^2 $. 

In Lemma \ref{oneparameter.lemma}
below we introduce sets ${\cal T}_{j , \epsilon} $ and ${\cal T}_{j, {\rm rem}}$ which we discuss
in Remark \ref{Tj.remark} following the lemma. The sup-script ``$\epsilon$" stands for the noise
term $\epsilon$. 
The subscript ``${\rm rem}$" stands
for ``remainder":\ under certain conditions terms with this subscript are of
smaller order than the other terms. 
 
\begin{lemma} \label{oneparameter.lemma} Let $\hat \gamma_j$ be obtained using the Lasso as described above. 
Let 
$$ {\cal T}_{j, \epsilon} := \{ | \hat c_j^T X^T \epsilon   |/n \le \lambda_{j, \epsilon} \hat \tau_j \} , \ 
{\cal T}_{j, {\rm rem} } := \{ 
 \lambda_j \| \hat \beta_{-j} - \beta_{-j}^0 \|_1   \le \lambda_{j,{\rm rem}}  \hat \tau_j\} .$$ 
On ${\cal T}_{j, \epsilon} \cap {\cal T}_{j, {\rm rem}}$ it holds that
$$| \hat \beta_j - \beta_j^0 | \le ( \lambda_{\epsilon} +  \lambda_{j, {\rm rem}})
{ \hat \tau_j / \tilde \tau_j^2}  + \lambda  \| \hat \theta_j \|_1 .
 $$
 The following result is moreover useful when $|\beta_j^0 |$ is small.
 If $\| \hat \gamma_j \|_1 < 1 $ and
% 
% , then on ${\cal T}_j\cap {\cal T}_{j, {\rm rem}}$
% $$  ( 1- \| \hat \gamma_j \|_1 ) | \hat \beta_j^0 | \le ( 1 + \| \hat \gamma_j \|_1 )
% |\beta_j^0 | +  { (\lambda_{\epsilon} + \lambda_{j, {\rm rem} } )^2 / (4 \lambda)}     $$
% and if
 $\lambda (1- \| \hat \gamma_j \|_1 ) > (\lambda_{j, \epsilon}  + \lambda_{j, {\rm rem}}) \hat \tau_j
  $, then on ${\cal T}_{j, \epsilon} \cap {\cal T}_{j, {\rm rem}}$
 $$\biggl  (\lambda (1- \| \hat \gamma_j \|_1 ) - (\lambda_{j, \epsilon} + \lambda_{j, {\rm rem}})\hat \tau_j \biggr )
  |\hat \beta_j | \le \biggl ( ( \lambda_{\epsilon}  + \lambda_{j, {\rm rem}} ) \hat \tau_j +
 \lambda \| \hat c_j \|_1 \biggr ) |\beta_j^0 | . $$
\end{lemma}

 \begin{remark} \label{Tj.remark} Note that if $\EE \epsilon=0 $ and $\EE \epsilon \epsilon^T = \sigma_{\epsilon}^2 I$,
 then 
${\rm var} ( \hat c_j^T X^T \epsilon/n )  = \sigma_{\epsilon}^2 \hat \tau_j^2/n $. Hence,
with $\lambda_{j, \epsilon}= {\mathcal O} (\sigma_{\epsilon} /\sqrt {n})$ large enough the set ${\cal T}_{j, \epsilon}$
will have large probability. The set ${\cal T}_{j,{\rm rem}}$ will e.g. have large probability
under compatibility conditions 
with $\lambda_{j, {\rm rem}} = {\mathcal  O} ( \lambda_j \lambda s / \hat \tau_j )$ where
$s$ is the number of non-zero $\beta_j^0$ or some
sparse approximation thereof, see e.g.\ \cite{brt09} and \cite{pbvdg11}, and see also
Corollary \ref{main.corollary}.
 \end{remark}

\begin{corollary} \label{oneparameter.corollary} Let us apply the square root Lasso
$$\hat \gamma_j := \arg \min_{\gamma_j \in \R^p} \biggl \{ 
\| X_j - X_{-j} \gamma_j\|_n + \lambda_0 \| \gamma_j \|_1 \biggr \} $$
where $\lambda_0$ is a ``universal"
tuning parameter (\cite{belloni2011square}). Then $\lambda_j = \lambda_0 \hat \tau_j $ and
so $\lambda_j \| \hat \beta_{-j} - \beta_{-j}^0\|_1 / \hat \tau_j =
\lambda_0 \| \hat \beta_{-j} - \beta_{-j}^0 \|_1 $ which can be bounded
by $\lambda_0 \| \hat \beta - \beta^0 \|_1 $. 
Define now ${\cal T}_{j , \epsilon}$ as in Lemma \ref{oneparameter.lemma} and 
${\cal T}_{\rm rem} := \{ \lambda_0 \| \hat \beta - \beta^0 \|_1 \le \lambda_{\rm rem} \} $.
We then have on ${\cal T}_{j, \epsilon}  \cap {\cal T}_{\rm rem} $
 $$\tilde \tau_j^2 | \hat \beta_j - \beta_j^0 |/ \hat \tau_j  \le \lambda_{j, \epsilon}  + \lambda_{\rm rem} +
 \lambda  \| \hat c_j \|_1/ \hat \tau_j    .$$
As noted in Remark \ref{Tj.remark} we typically can take 
$\lambda_{\rm rem} \asymp
\lambda_0 \lambda s $ and moreover typically
$\lambda_0 = {\mathcal O} (\sqrt {\log p / n } )$ and
$\lambda = {\mathcal O}_{\PP}
(\sigma_{\epsilon} \sqrt { \log p / n})  $.
Hence, under suitable conditions on the sparseness $s$ the remainder term $\lambda_{\rm rem}$ is negligible.
Note that the first term $\lambda_{j, \epsilon}$ needs to be chosen be proportional on the standard deviation $\sigma_{\epsilon}$ of the noise $\epsilon$
to ensure that ${\cal T}_{j, \epsilon}$ has large probability.  The tuning $\lambda$ has to be chosen proportional
to (an estimator of) $\sigma_{\epsilon}$ as well. Thus, all three terms scale with $\sigma_{\epsilon}$.

\end{corollary}

{\bf Standard asymptotic scenario I} 
For a better reading of the bounds we present
them in an asymptotic formulation.  In all asymptotic formulations in this section we
assume 
$\EE \epsilon =0 $ and $\EE \epsilon \epsilon^T = \sigma_{\epsilon}^2 I$
where $\sigma_{\epsilon} = {\mathcal O} (1)$. 
Suppose that $\hat \tau_j^2 = {\mathcal O} (1)$ (which is true
of the columns of $X$ are normalized so that $\| X_j \|_n=1$), $1/ \hat \tau_j^2 = {\mathcal O} (1)$
(a restricted eigenvalue or compatibility condition), 
 $\lambda =  {\mathcal O}_{\PP}
 (\sqrt {\log p / n } )$ (the standard choice), that with $\lambda_{j, \epsilon}= {\mathcal O}_{\PP}
  (n^{-1/2})$ and $\lambda_{j, {\rm rem} } =
 {\mathcal O} (n^{-1/2})$
 the probability
 of ${\cal T}_{j, \epsilon} \cap {\cal T}_{j, {\rm rem}} $ goes to one
(this follows from the moment conditions on the noise and when
$\lambda_0 = {\mathcal O}
 ( \sqrt {\log p / n} )$ from $s = {\mathcal O} ( \sqrt n / \log p )$ and a compatibility condition).
 Then Lemma \ref{oneparameter.lemma} gives that
 $$| \hat \beta_j - \beta_j^0 | =
  {\mathcal O}_{\PP} (n^{-1/2}) +
 \lambda \| \hat \theta_j \|_1 $$
 and if $\| \hat \gamma_j \|_1 < 1 $
 $$| \hat \beta_j | \le { \| \hat c_j \|_1 + o_{\PP}(1)\over 
 1- \| \hat \gamma_j \|_1 + o_{\PP}( 1) } | \beta_j^0 | . $$

\subsection{A bound for the (weighted) $\ell_{\infty}$-norm}\label{infinity.section}

The above can be applied for each $j$. Let $\hat C:= ( \hat c_1 , \ldots , \hat c_p ) $, 
$\hat \Theta := ( \hat \theta_1 , \ldots , \hat \theta_p )$,
and $\hat T:= {\rm diag} ( \hat \tau_1 , \ldots , \hat \tau_p )$, $\tilde T := {\rm diag} ( \tilde \tau_1 , \ldots ,
\tilde \tau_p) $. 
Define for a matrix $A$ the $\ell_1$-operator norm
$$  \vvert A \vvert_1 : = \max_{j} \sum_{k} | A_{k,j} | .$$

We only present the result when using the square root Lasso because of
its elegant scaling.
We obtain a bound for the $\ell_{\infty}$-estimation error in terms of the
$\ell_1$-estimation error. The weighted $\ell_{\infty}$-norm will be of interest
when the residual variances $\hat \tau_j^2$ are not balanced for
the various values of $j$. Note also that we may use the bound
$\hat \tau_j/ \tilde \tau_j \le 1$ for all $j$. 

\begin{lemma} \label{all.lemma} Define
$${\cal T}_{\rm all} := \{ \| \epsilon^T X \hat C \hat T^{-1}  \|_{\infty} \le \lambda_{\epsilon}  \} . $$
When using the square root Lasso as in 
Corollary \ref{oneparameter.corollary}, we have on ${\cal T}_{\rm all}$
$$ \| \hat \beta - \beta^0 \|_{\infty} \le  ( \lambda_{\epsilon}   + \lambda_0 \| \hat \beta - \beta^0 \|_1) \| 
 \tilde T^{-2}\hat T  \|_{\infty}  +
\lambda \vvert \hat \Theta \vvert_1   $$
and
$$\|  \hat T^{-1} \tilde T^2 ( \hat \beta - \beta^0 ) \|_{\infty} \le  \lambda_{\epsilon}  +\lambda_0 \| \hat \beta - \beta^0 \|_1 +
\lambda \vvert  \hat C \hat T^{-1} \vvert_1  . $$
\end{lemma}

{\bf Standard asymptotic scenario II} Assume that $\| \hat T^{-1} \|_{\infty} =
{\mathcal O} (1)$ and that for 
$\lambda_{\epsilon} = {\mathcal O} (\sqrt {\log p /n})$ the probability of ${\cal T}_{\rm all}$
goes to one and for $\lambda={\mathcal O}_{\PP} (\sqrt {\log p / n})$ 
it holds that 
$\lambda_0 \| \hat \beta - \beta^0 \|_1 = {\mathcal O} _{\PP}(\sqrt {\log p / n  })$.
Then we get
$$\| \hat \beta - \beta^0 \|_{\infty} = {\mathcal O}_{\PP} ( \sqrt {\log p / n} )  +
\lambda \vvert \hat \Theta \vvert_1 . $$

\subsection{De-sparsifying the Lasso and re-sparsifying}\label{de-sparsity.section}

The de-sparsified Lasso is defined in \cite{vdgetal13} as
$$\hat b_j = \hat \beta_j + \hat \theta_j^T X^T (Y- X \hat \beta) /n , \ j=1 , \ldots , p $$
or in matrix notation
$$ \hat b = \hat \beta + \hat \Theta^T X^T ( Y- X \hat \beta )/n . $$

The de-sparsified Lasso removes the
term involving the $\ell_1$-norm of the worst possible sub-direction.
In that sense, it removes the bias due to the $\ell_1$-penalty
on this sub-direction. We prove this
for completeness, but the result is as in \cite{vdgetal13} where also the details concerning the resulting (asymptotic)
normality of the de-sparsified Lasso are presented.

\begin{lemma} \label{de-sparsify.lemma} Let $\hat \Theta $ be obtained using the square root Lasso.
Then
$$\hat T^{-1} \tilde T^2  ( \hat b - \beta^0) = \hat T^{-1} \tilde T^2 \hat \Theta^T X^T \epsilon/n +
{\rm rem} , $$
where $\| {\rm rem}\|_{\infty}  \le  \lambda_0 \| \hat \beta - \beta_0 \|_1$.
Hence for a fixed $j$, on 
$${\cal T}_{j, \epsilon} \cap {\cal T}_{\rm rem} := \{ | \hat c_j^T X^T \epsilon | /n \le \lambda_{j, \epsilon } \hat \tau_j , \ \lambda_0 \| \hat \beta - \beta^0 \|_1 \le \lambda_{\rm rem} \} ,$$
we have
$$  \tilde \tau_j^2 | \hat b_j - \beta_j^0 |/  \hat \tau_j  \le \lambda_{j, \epsilon} +
\lambda_{\rm rem}  $$
and on 
$${\cal T}_{\rm all} \cap {\cal T}_{\rm rem}  :=  \{ \| \hat T^{-1} \hat C^T X^T \epsilon    \|_{\infty} /n \le \lambda_{\epsilon}  , \ \lambda_0 \| \hat \beta - \beta^0 \|_1 \le \lambda_{\rm rem} \} $$
 it holds
that
$$\| \hat T^{-1} \tilde T^2   ( \hat b - \beta^0 ) \|_{\infty} \le  \lambda_{\epsilon}  +\lambda_{\rm rem} . $$\end{lemma}

{\bf Re-sparsifying}
One may want to re-sparsify the de-sparsified Lasso $\hat b$ using some threshold
$\lambda_{\rm sparse}$ giving the estimator
$$\hat b_{j, {\rm sparse}} := \hat b_j {\rm l} \{ | \hat b_j | > \lambda_{\rm sparse} \hat \tau_j /
\tilde \tau_j^2 \} , \ j=1 , \ldots , p . $$
As we see from Lemma \ref{de-sparsify.lemma} this new estimator can improve the $\ell_{\infty}$-bounds of the Lasso and
has under sparsity conditions $\ell_q$-bounds similar to the Lasso ($1 \le q < \infty$).

{\bf Standard asymptotic scenario} Under Scenario I we have
$$| \hat b_j - \beta_j^0 | = {\mathcal O}_{\PP} (n^{-1/2}) , $$
and under Scenario II 
$$\| \hat b - \beta^0 \|_{\infty} = {\mathcal O}_{\PP} (\sqrt {\log p / n } )  . $$
This can be used for exact recovery of the relevant active set where the coefficients 
are sufficiently larger than $\sqrt {\log p /n }$ in absolute value. In other words, it leads
to exact recovery without assuming an irrepresentable condition.
For the re-sparsified estimator we also take $\lambda_{\rm sparse} \asymp 
\sqrt {\log p /n }$. Then under Scenario II
$$ \| \hat b_{\rm sparse} - \beta^0 \|_{\infty} = {\mathcal O}_{\PP} (\sqrt {\log p / n } ) .$$
Assuming as an example that for some 
fixed $r \in (0,1)$ it holds that $\sum_{j=1}^p |\beta_j^0 |^r = {\mathcal O}(1)$, then
we obtain under Scenario II
$$\| \hat b_{\rm sparse}- \beta^0 \|_q = {\mathcal O}_{\PP}
(s^{1 \over q}  \sqrt {\log p / n } ) , \ 1 \le q \le \infty $$ 
where $s = (\sqrt {n / \log p } )^r $.

\subsection{Using the Dantzig selector} \label{Dantzig.section}
Instead of defining approximate projections with help of the Lasso, we may also
use the Dantzig selector, which gives a new definition of $\hat \gamma_j$:
$$ \hat \gamma_j := \arg \min_{\gamma_j \in \R^{p-1} } \| \gamma_j \|_1 \ {\rm s.t.} \
\| X_{j}^T ( X_j - X_{-j} \gamma_j ) \|_{\infty} / n \le \lambda_j . $$

The result of 
Lemma \ref{Dantzig.lemma} below is quite similar to the one of Corollary \ref{oneparameter.corollary} .

\begin{lemma} \label{Dantzig.lemma} Let $\hat \gamma_j$ be obtained using the Dantzig selector. 
Assume $\lambda_j \| \hat \gamma_j \|_1 < \hat \tau_j^2$. On 
$${\cal T}_{j, \epsilon}
\cap {\cal T}_{j, {\rm rem} } 
:= \{ | \hat c_j^T X^T \epsilon   |/n \le \lambda_{j, \epsilon} \hat \tau_j , \ 
 \lambda_j \| \hat \beta_{-j} - \beta_{-j}^0 \|_1/ \hat \tau_j \le \lambda_{j, {\rm rem}} 
\}$$
we have
$$\biggl  [ \hat \tau_j^2 - \lambda_j \| \hat \gamma_j \|_1  \biggr ] { | \hat \beta_j - \beta_j^0 |
 \over \hat \tau_j } \le \lambda_{j, \epsilon} + \lambda_{j, {\rm rem}} +
 \lambda
 \| \hat c_j \|_1   /\hat \tau_j  .$$
 \end{lemma}

\subsection{A bound for a group of variables}\label{group.section}
Recall the $\ell_1$ operator norm
$\vvert A \vvert_1 := \max_j \sum_{k} | a_{k,j} | $ of a matrix $A$.
We define its $\ell_1$ norm as
$\| A \|_1 := \sum_{j} \sum_{k} |a_{k,j} | $.
Let $J \subset \{ 1 , \ldots , p \}$ be a given subset of the variables.
We let $X_J := \{ X_j \}_{j \in J}$ and $X_{-J} := \{ X_j \}_{j \notin J } = X_{J^c} $.
Moreover, we write $\beta_J := \{\beta_j \}_{j \in J} \in \R^{|J|} $ and we use
the same notation for $\{ \beta_j {\rm l} \{ j \in J \} \}_{j=1}^p  \in \R^p $ and likewise
for $\beta_{-J} = \beta_{J^c} $. 
We let 
$$ \hat \Gamma_J := \arg \min_{\Gamma_J }\biggl \{ 
{\rm trace} ( X_J - X_{-J} \Gamma_J )^T ( X_J - X_{-J} \Gamma_J ) + \lambda_{J,0} \|  \Gamma_J \|_1
\biggr \}  . $$
In other words, each column $\hat \gamma_{J,j}$ of $\hat \Gamma_J$ is obtained
by performing a Lasso of $X_j$ on $X_{-J}$ with tuning parameter $\lambda_{J,0}$.
We define $\hat C_J $ by $X \hat C_J = X_J -
X_{-J} \hat \Gamma_J $ so that $\vvert \hat C_J \vvert_1 = 1 + \vvert
\hat \Gamma_J \vvert_1 $. 
We introduce the smallest eigenvalue
$$\hat \phi_J^2 := \min \{ \| X \hat C_J \beta_J \|_n^2: \ \| \beta_{J} \|_2 = 1 \}  .$$
The compatibility constant is
$$ \hat \phi^2 ( L, S) := \min \{|S|  \| X \beta \|_n^2 : \ \| \beta_{-S } \|_1 \le L , \| \beta_S \|_1 = 1 \} , $$
see \cite{vandeGeer:07a} or \cite{pbvdg11}. 

Let us make three remarks. Firstly, 
we note that  $\vvert \hat \Gamma_J \vvert_1$
is generally not directly comparable to
$\max_{j \in J} \|\hat \gamma_j \|_1$ ($\hat \gamma_j$, $j =1 , \ldots, p$, defined in 
Subsection \ref{single.section}). Secondly, 
in view of the scaling in  Part I of Theorem \ref{main.theorem} below, a matrix version of the square root Lasso would be to take
$\lambda_{J,0} = \lambda_0 \hat \phi_J $. And thirdly, 
it is easy to see that
$$\hat \phi_J^2 \le \hat \phi^2 (L , J\cap S ) , \forall \ L, S .$$
Theorem \ref{main.theorem} below may be applied to
general and hence also relatively small sets $J $. In that case one may want to replace
$\hat \phi ( L, J \cap S)$ by 
by $\hat \phi_J$ so that the constant $L$ defined in the theorem  no longer plays any role.
(Moreover, we let $\hat \phi(\infty,S)= \hat \phi_J$).

In Part I of the theorem we establish bounds under relatively weak (see
(\ref{for-desparsifying.equation})) or no (see (\ref{simpler.equation})) conditions
on the (approximate) worst possible  sub direction $\hat \Gamma_J$. Part II
of the theorem assumes $\vvert \hat \Gamma_J \vvert_1$ is sufficiently small
and presents oracle bounds.

\begin{theorem}  \label{main.theorem} 
Let $\hat \Gamma_J$ be obtained using the Lasso.\\
{\bf Part I} Define
$$\bar {\cal T}_{J, \epsilon} := \{ | \beta_J^T \hat C_J^T X^T \epsilon |/n \le \bar \lambda_{J, \epsilon} \sqrt {|J|} 
\| X \hat C_J \beta_J \|_n , \ \forall \ \beta_J \in \R^{|J|} \} . $$
 If
$\lambda_{J,0} |J| \vvert \hat \Gamma_J \vvert_1 < \hat \phi_J^2 $,
then we have on $\bar {\cal T}_{J, \epsilon}$
\begin{equation}\label{for-desparsifying.equation}
\biggl ( 1- { \lambda_{J,0} |J| \vvert \hat \Gamma_J \vvert_1 \over
\hat \phi_J^2 }  \biggr ) \| X \hat C_J ( \hat \beta_J - \beta_J^0 ) \|_n \end{equation}
$$ \le {\sqrt{ |J| } \over \hat \phi_J } \biggl (\bar \lambda_{ J,\epsilon} \hat \phi_J +  \lambda_{J,0} \| \hat \beta_{-J}-
\beta_{-J}^0 \|_1 
+ \lambda \vvert \hat C_J \vvert_1  \biggr ).$$
Furthermore on $\bar {\cal T}_{J, \epsilon}$
\begin{equation}\label{simpler.equation}
 \| X \hat C_J ( \hat \beta_J - \beta_J^0 ) \|_n \le
{\sqrt {|J|} \over \hat \phi_J } \biggl ( 
\bar \lambda_{J, \epsilon} \hat \phi_J +  (\lambda + \lambda_{J, 0}  \| \hat \beta- \beta^0 \|_1 )
\vvert \hat C_J \vvert_1 \biggr  ) 
\end{equation}
which holds without assuming some bound for $\vvert \hat \Gamma_J \vvert_1$.\\
{\bf Part II} 
Define
$${\cal T}_{J, \epsilon} := \{\| \hat C_J^T X^T \epsilon\|_{\infty} /n \le \lambda_{ J, \epsilon} \} $$
and
$${\cal T}_{J, {\rm rem}} := \{ \lambda_{J,0} \| \hat \beta - \beta^0 \|_1 \le \lambda_{J, {\rm rem}} \} . $$
Assume $\vvert \hat \Gamma_J \vvert_1  < 1 $ and that in fact for some $\lambda_1\ge 0$
$$\lambda (1- \vvert \hat \Gamma \vvert_1 ) \ge  \lambda_{J,\epsilon} + \lambda_{J, {\rm rem}} \vvert \hat C_J \vvert_1 + \lambda_1.$$
Let
$$ L:= {  \lambda_{J, \epsilon} +( \lambda+ \lambda_{J,{\rm rem}}   )  \vvert \hat C_J  \vvert_1 + \lambda \vvert \hat \Gamma_J \vvert_1 + \lambda_1 
\over \lambda ( 1- \vvert \hat \Gamma_J \vvert_1 ) - ( \lambda_{J, \epsilon} + \lambda_{J,{\rm rem}}  \vvert \hat
C_J \vvert_1 ) - \lambda_1 }  . $$
Let $\beta \in \R^p$ be arbitrary and let $S:= \{ j :\ \beta_j \not= 0 \}$ be its active set. 
Then on ${\cal T}_{J, \epsilon} \cap {\cal T}_{J, {\rm rem}}$,
$$\| X \hat C_J ( \hat \beta_J - \beta_J^0 ) \|_n^2 +2 \biggl [ \biggl (\lambda (1- \vvert \hat \Gamma_J \vvert_1 )  -(\lambda_{ J, \epsilon}  +  \lambda_{J,{\rm rem}}   \vvert \hat C_J \vvert_1 ) - \lambda_1 \biggr) \wedge
\lambda_1 \biggr ]  \| \hat \beta_{J }- \beta_J\|_1  $$ $$
  \le {{4 \lambda^2  | J \cap S | } \over \hat \phi ^2 (L, J \cap S) }
 + \| X \hat C_J ( \beta_J - \beta_J^0) \|_n^2  . $$
\end{theorem}

\begin{remark} Although (\ref{simpler.equation}) is a rougher bound than
 (\ref{for-desparsifying.equation})
(assuming the condition on $\vvert \hat \Gamma_J \vvert_1$)
its simplicity makes it preferable. We have stated 
(\ref{for-desparsifying.equation}) because
for de-sparsifying we actually need this more refined result
(see Lemma \ref{desparsify-group.lemma}).

\end{remark} 

\begin{corollary} \label{main.corollary} $ $\\
{\bf Sharp oracle inequality} 
If we take $J= \{ 1 , \ldots , p \} $ (and $\lambda_{J,0}=0$)
we recover from Part II of Theorem \ref{main.theorem} the sharp oracle inequality: on ${\cal T}_{\epsilon}:=
\{ \| \epsilon^T X \|_{\infty} /n\le \lambda_{\epsilon}\} $ with $\lambda_{\epsilon} < \lambda$:
$$\| X ( \hat \beta - \beta^0 ) \|_n^2 + 2 \biggl ( ( \lambda -\lambda_{\epsilon}- \lambda_1
) \wedge \lambda_1  \biggr )\| \hat \beta - \beta \|_1 
\le { 4 \lambda^2 { |  S_0 | } \over \hat \phi^2 (L,  S) }
 + \| X ( \beta - \beta^0 ) \|_n^2  $$
with $L= (\lambda + \lambda_{\epsilon}+ \lambda_1  ) /( \lambda - \lambda_{\epsilon} - \lambda_1 )$.
This corresponds to results in \cite{koltchinskii2011nuclear} or \cite{vdG2014}.\\
{\bf Small coefficients}
If we take $J = S^c$ where $S$ is such that
${\rm rank} (X_S )= |S|$, and we choose $\lambda_{J,0} =0$ we find from Part II of Theorem \ref{main.theorem}: on ${\cal T}_{\epsilon}$: under the irrepresentable condition on the set $S^c$
$$ \vvert \hat \Gamma_{S^c} \vvert_1 \le ( \lambda - \lambda_{\epsilon}  ) / ( \lambda + \lambda_{\epsilon})$$
we have
\begin{equation}\label{betasmall.equation}
\biggl (\lambda - \lambda_{\epsilon} - \vvert \hat \Gamma_J \vvert_1 (\lambda + \lambda_{\epsilon} )
\biggr )  \|\hat \beta_{S^c} \|_1 \le \| X\hat C_{S^c}   \beta_{S^c}^0 \|_n^2 .
 \end{equation}
 (We took $2 \lambda_1= \lambda- \lambda_{\epsilon}- \vvert \hat \Gamma_J \vvert_1(\lambda +\lambda_{\epsilon} ) $.)
 This bound generalizes the bound for small
 values of $| \beta_j^0 |$ presented in Lemma \ref{oneparameter.lemma}. 
To bound  $\|X \hat C_{S^c} \beta_{S^c}^0  \|_n^2$ one may want to use 
$$\| X\hat C_{S^c}  \beta_{S^c}^0 \|_n^2 \le  \| \beta_{S^c}^0 \|_2^2
\| \hat C_{S^c} X^T X \hat C_{S^c}\|_{2,2}^2 / n
 $$
where $\| A \|_{2,2}^2$ is the largest eigenvalue of a positive semi-definite matrix $A$.
Moreover
\begin{equation}\label{maxeigenvalue.equation}
\| \hat C_{S^c} X^T X \hat C_{S^c}\|_{2,2}^2 / n \le \| X_J^T X_J \|_{2,2}^2/n + 2 \lambda_{J,0} |J|  . 
\end{equation}
Alternatively (and for a comparison with Theorem \ref{varselection.theorem} below)
one may use the bound
$$\| X \hat C_{S^c} \beta_{S^c}^0  \|_n^2 \le
 \| \beta_{S^c}^0 \|_1 \| \| X \hat C_{S^c} \beta_{S^c}^0 \|_{\infty}/n \le
 \| \beta_{S^c}^0 \|_1 \| \beta_{S^c}^0 \|_{\infty}  \vvert \hat C_{S^c}^T X^T X \hat C_{S^c}   \vvert_1 /n .$$
 {\bf Single coefficients} 
If we take $J = \{ j \}$, $\lambda_{J,0} := \lambda_j$ and use
the first result (\ref{for-desparsifying.equation}) of Part I of Theorem \ref{main.theorem} we find on ${\cal T}_{j, \epsilon} \cap
{\cal T}_{j, {\rm rem}} := \{| \hat c_j^T X^T \epsilon |/n \le \lambda_{\epsilon}
\hat \tau_j , \ \lambda_j \| \hat \beta_{-j}- \beta_{-j}^0 \|_1 /\hat \tau_j \le \lambda_{j, {\rm rem}} \} $
$$( \hat \tau_j^2 - \lambda_j  \| \hat \gamma_j \|_1 )  | \hat \beta_j - \beta_j^0 |/ \hat \tau_j  \le  
  \lambda_{\epsilon} +\lambda_{j, {\rm rem}} + \lambda \| \hat c_j \|_1 / \hat \tau_j  $$
which is similar to Lemma \ref{Dantzig.lemma}.  With a more refined handling of 
the cross terms in the matrix computations
one can of course also  recover the result of Lemma \ref{oneparameter.lemma}.
\end{corollary}

{\bf Standard asymptotic scenario III}\\
The constant $\bar \lambda_{J, \epsilon}$ in Part I of Theorem
\ref{main.theorem} can be taken of order $1/ \sqrt n $.
If the design is normalized (${\rm diag} (X^T X)/n = I$) 
the constant $\lambda_{J, \epsilon}$ in Part II can generally be taken of order
$1/\sqrt n$ when $J$ is finite. Else we take $\lambda_{J, \epsilon} =
\lambda_{\epsilon} \vvert \hat C_J \vvert_1 $ with $\lambda_{\epsilon}=
{\mathcal O} (\sqrt {\log p / n } )$. Note that $\hat C_J^T X^T \epsilon/n$ has covariance
matrix $\hat C_J^T X^T X \hat C_J/n$, whose maximal eigenvalue can be bounded
as in (\ref{maxeigenvalue.equation}). 
 We suppose that the probability of ${\cal T}_{\epsilon} := \{ \| \epsilon^T X \|_{\infty}/n  \le \lambda_{\epsilon}\} $
 goes to one for $\lambda_{\epsilon}={\mathcal O}_{\PP} (\sqrt {\log p / n})$
 suitably chosen. We assume $\lambda =  {\mathcal O}_{\PP} (\sqrt {\log p / n})$ and $\lambda_{J, {\rm rem}} =
{\mathcal O}_{\PP} (\sqrt {\log p/ n})$, that $1/ \hat \phi_J = {\mathcal O}(1)$
and that the largest eigenvalue of $X^T X  /n $ is ${\mathcal O} (1)$.
Then we find from Part II of Theorem \ref{main.theorem}, under the assumption
$\vvert \hat \Gamma_J \vvert_1 \le (\lambda - \lambda_{\epsilon}) / (\lambda + \lambda_{\epsilon}) +
o(1)$, that
$$ \| \hat \beta_J - \beta_J^0 \|_2^2 = {\mathcal O}_{\PP} \biggl ( \log p | J \cap S| /n + \| \beta_J -
\beta_J^0 \|_2^2 \biggr ) , $$
where for suitable $\lambda_{\rm thres} = {\mathcal O} ( \sqrt {\log p / n } ) $
the vector $\beta$ is the sparse approximation
$$\beta_{j} := \beta_j^0 \{ | \beta_j^0 | > \lambda_{\rm thres} \} , \ \forall \ j  , 
$$ and  $S := \{ j :\ |\beta_j^0 | >\lambda_{\rm thres} \} $ is its active set.

\subsection{Variable selection in the
approximately sparse case}\label{selection.section}

The irrepresentable condition on the set of inactive variables is commonly used to show that the 
$\ell_1$-penalized 
estimator $\hat \beta$ has no false positives.
Inequality (\ref{betasmall.equation}) in Corollary \ref{main.corollary} shows that 
under irrepresentable conditions on a set $S^c$ where the coefficients
of $\beta_j^0$, $j \in S^c$, are small, the estimated coefficients $|\hat \beta_j|$, $j \in S^c$
will be small as well. We now show
that under stronger bounds for $\vvert \hat \Gamma_{S^c} \vvert_1$ actually  $\hat \beta_j$
will be zero for $j \in S^c$.  This result thus extends the situation to the approximately sparse
case where there may be many non-zero but small coefficients. It can be a step
towards local uniformity and away from super-efficiency and may be useful
for building confidence intervals.
Note that $X \hat C_{S^c} \beta_{S^c}^0 $ is the
part of $X \beta_{S^c}^0$ left over after projecting it on 
$ X_S $.

\begin{theorem} \label{varselection.theorem} Let $\hat \beta$ be the unique solution of the KKT conditions
(\ref{KKT.equation}). Let
${\cal T}_{\epsilon} := \{ \| X^T \epsilon \|_{\infty} / n \le \lambda_{\epsilon} \} $.
Consider some set $S \subset \{ 1 , \ldots , p \}$ with
${\rm rank} (X_S) = |S|$ and define $\hat \Gamma_{S^c} := (X_S^T X_S)^{-1} X_S^T X_{S^c} $
and $X \hat C_{S^c} = X_{S^c} - X_S \hat \Gamma_{S^c} $. Suppose 
$$ (\lambda + \lambda_{\epsilon} )\vvert \hat \Gamma_{S^c} \vvert_1  + \| \hat C_{S^c}^T X^T X \hat C_{S^c} \beta_{S^c}^0 \|_{\infty}/n \le   \lambda - \lambda_{\epsilon } 
 . $$
Then on ${\cal T}_{\epsilon}$ we have $\hat \beta_j = 0 $ for all $j \notin S$. 

\end{theorem}

\subsection{De-sparsifying a group of variables}\label{de-group.section}

We define the group de-sparsified estimator
$$ \hat b_J := \hat \beta_J + ( \hat C_J^T X^T X_J )^{-1} \hat C_J^T X^T (Y - X \hat \beta )  $$
assuming the above used matrix inverse exists. 

Again we shall need the smallest
eigenvalue
$$\hat \phi_J^2 := \min \{ \| X \hat C_J \hat \beta_J \|_n^2:\ \| \beta_J \|_2 = 1 \}  $$
which is bounded from below
by the compatibility constant $\hat \phi^2 ( \vvert \hat \Gamma \vvert_1 , J)$.

\begin{lemma}\label{desparsify-group.lemma} 
Let 
$$\bar {\cal T}_{J, \epsilon} := \{ | \beta_J^T \hat C_J^T X^T \epsilon|  /n \le \bar \lambda_{J, \epsilon}
\sqrt {|J|} 
\| X \hat C_J \beta_J \|_n , \forall \ \beta_J \in \R^{|J|}  \} .$$
If $\lambda_{J,0} | J| \vvert \hat \Gamma_J \vvert_1 / \hat \phi_J < 1 $,
then on $\bar {\cal T}_{J, \epsilon}$
$$ \biggl ( 1- { \lambda_{J,0} |J| \vvert \hat \Gamma_J \vvert_1 \over
\hat \phi_J } \biggr) \| X\hat C_J ( \hat b_J - \beta_J^0) \|_n \le 
{ \sqrt {|J|} \over \hat \phi_J } \biggl ( \bar \lambda_{J, \epsilon} \hat \phi_J +
\lambda_{J,0} \| \hat \beta_{-J} - \beta_{-J}^0 \|_1 \biggr ) . $$
Let furthermore
$$ {\cal T}_{J, \epsilon}^* := \{ \|  \hat C_J^T X^T \epsilon \|_2 /\sqrt n \le
\sqrt {|J|}  \lambda_{J, \epsilon}^* \} .$$
Then on ${\cal T}_{J, \epsilon}^* $,
$$\| \hat C_J^T X^T X_J ( \hat b_J - \beta_J^0 ) \|_2 /\sqrt n \le
{ \sqrt {|J|} \over \hat \phi_J}  \biggl ( \lambda_{J, \epsilon}^* \hat \phi_J + \lambda_{J,0} \| \hat \beta_{-J} - \beta_{-J}^0 \|_1 \biggr  ) 
.$$
\end{lemma}

{\bf Standard asymptotic scenario IV} Suppose that $1/ \hat \phi_J = {\mathcal O} (1)$ and that for a suitable $\lambda =
{\mathcal O}_{\PP} ( \sqrt {\log p / n })$ it holds that
$\lambda_{J,0} \| \hat \beta - \beta^0 \|_1 = {\mathcal O}_{\PP}  ( n^{-1/2})$ 
(recall this typically holds if approximate sparsity is of
order $\sqrt n / \log p $) and that for suitable
$\lambda_{J, \epsilon}^*= {\mathcal O}_{\PP}  (n^{-1/2})$
the probability of $ {\cal T}_{J, \epsilon}^* $ tends to one. Then
$$ \|  \hat C_J^T X^T X_J ( \hat b_J - \beta_J^0 ) \|_2 /\sqrt n = {\mathcal O}_{\PP}
(\sqrt {|J|/n}  ) . $$

If we assume moreover for suitable
$\bar \lambda_{J, \epsilon}= {\mathcal O}_{\PP}  (n^{-1/2})$
the probability of $ \bar {\cal T}_{J, \epsilon} $ tends to one  and
that $\vvert \hat \Gamma_J \vvert_1 \sqrt {|J|} \sqrt {\log p / n} / \hat \phi_J  = {\mathcal O} (1)$
is suitably small, 
then
$$ \| \hat b_J - \beta_J^0 \|_2 = {\mathcal O}_{\PP} ( \sqrt { |J| / n } ) . $$

\section{Random design}\label{random.section}
We study the linear model (\ref{linearmodel.equation}) but now with $X$ a random
matrix with i.i.d.\ rows with distribution $P$. Again $\hat \beta$ is the Lasso
estimator defined in (\ref{Lasso.equation}). 
We let $X_{-J} \Gamma_J$ be the projection of $X_J$ on $X_{-J}$ in $L_2 (P)$
and $X C_J := X_J - X_{-J} \Gamma_J $.
Then by the same arguments as in Theorem \ref{main.theorem} we find 
Theorem \ref{mainrandom.theorem} below. 
To avoid digressions we omit the counterpart of the first inequality 
(\ref{for-desparsifying.equation}) in Theorem \ref{main.theorem} as 
it has no direct counterpart for de-sparsifying 
since we now use
theoretical projections which are usually unknown (as $P$ is usually unknown).

With some abuse of notation, we let $\hat \phi_J ^2$ now be the smallest
eigenvalue of the matrix $C_J^T X^T X C_J / n $. We note that the formulation in
Theorem \ref{mainrandom.theorem} is again in terms of the (now random) norm
$\| \cdot \|_n$ and (now random) compatibility constants $\hat \phi(L,S)$.
In our results of the next section for general loss functions we use an alternative approach,
leading in the case of random design to a formulation where the
difference between the $L_2(P)$ norm and the empirical norm $\| \cdot \|_n$
ends up in the remainder term (see Subsection \ref{rem-linear.section}).

\begin{theorem}  \label{mainrandom.theorem} 
Define
$$ \bar {\cal T}_{J,\epsilon} := \{ |\beta_J^T C_J^T X^T \epsilon |/n \le\bar \lambda_{J, \epsilon} \sqrt {|J|} 
\| X C_J \beta \|_n , \ \forall \ \beta_J \in \R^{|J|} \} , $$
$${\cal T}_{J, \epsilon} := \{\| C_J^T X^T \epsilon\|_{\infty} /n \le \lambda_{ J, \epsilon}\}, $$
$$
{\cal T}_{J,{\rm rem}} :=  \{  \| X_{-J}^T X C_J\|_{\infty}  \| \hat \beta - \beta^0 \|_1 / n \le \lambda_{J,{\rm rem} } \} . $$
We have on $\bar {\cal T}_{J, \epsilon} \cap {\cal T}_{J, {\rm rem}} $
 $$ \| X C_J ( \hat \beta_J - \beta_J^0 ) \|_n \le 
{ \sqrt {J} \over \hat \phi_J } \biggl [ \bar  \lambda_{J, \epsilon} \hat \phi_J +( \lambda
 +  \lambda_{J,{\rm rem}}) \vvert C_J \vvert_1 \biggr ] 
.$$
Assume $\vvert \Gamma_J \vvert_1  < 1 $ and that in fact that for some $\lambda_1 \ge 0$
$$\lambda (1- \vvert \Gamma \vvert_1 - \lambda_1 ) >  \lambda_{J, \epsilon} + \lambda_{J,{\rm rem}}  \vvert  C_J \vvert_1 + \lambda_1.$$
Let
$$ L:= {  \lambda_{J, \epsilon} + (\lambda + \lambda_{J,{\rm rem}}   )  \vvert C_J \vvert_1  + \lambda \vvert \Gamma_J \vvert_1 + \lambda_1 
\over \lambda ( 1- \vvert  \Gamma_J \vvert_1 ) - ( \lambda_{J, \epsilon} + \lambda_{J,{\rm rem}}   \vvert 
C_J \vvert_1 ) - \lambda_1 }  . $$
Then on ${\cal T}_{J, \epsilon} \cap {\cal T}_{J, {\rm rem}}$,
$$\| X C_J ( \hat \beta_J - \beta_J^0 ) \|_n^2 + 2 
\biggl [ \biggl (\lambda (1- \vvert  \Gamma_J \vvert_1 )  -(\lambda_{ J, \epsilon}  +  \lambda_{J,{\rm rem}}   \vvert  C_J \vvert_1 ) - \lambda_1 \biggr ) \wedge
\lambda_1 \biggr ]  \| \hat \beta_{J }- \beta_J\|_1 
$$ $$
\le 
 { 4 \lambda^2 | J \cap S |  \over \hat \phi^2 (L, J \cap S) } + \| X  C_J ( \beta_J - \beta_J^0 ) \|_n^2 .
   $$
 \end{theorem}

\section{General loss functions} \label{general.section}
Let $X_1 , \ldots , X_n$ be independent observations in some observation space
${\cal X}$. For $d\in \Nat$ and a function $f: \ {\cal X} \rightarrow \R^d$ we use the notation
$$P_n f := {1 \over n} \sum_{i=1}^n f(X_i) , \ Pf := \EE P_n f . $$
Consider a convex subset $\Theta$ of $\R^p$ and a loss
function $\rho_{\theta}: \ {\cal X} \rightarrow \R$ with
derivative $\dot \rho_{\theta}:= \partial \rho_{\theta} / \partial \theta $,
$\theta \in \Theta$.
We examine the $\ell_1$-penalized M-estimator
$$\hat \theta := \arg \min_{\theta \in \Theta }
\biggl \{ P_n \rho_{\theta} + \lambda \| \theta \|_{1, {\rm off} } \biggr \} , $$
where $\| \theta \|_{1, \rm off} = \sum_{j \notin {\cal J}} | \theta_j | $
is the $\ell_1$-penalty on the parameters $\{ \theta_j \}_{j \notin {\cal J}}$
and ${\cal J} \subset \{1 , \ldots , p \}$ is a fixed set.
The subscript ``${\rm off}$" refers to the set of parameters $\{ \theta_j \}_{j \notin {\cal J} }$
which are susceptible to being turned off (i.e. set to zero). The set ${\cal J}$ contains
the indexes of parameters that are not penalized. 
We assume that 
 $\hat \theta$ is the unique solution of the KKT conditions
\begin{equation}\label{KKT2.equation}
P_n \dot \rho_{\hat \theta} + \lambda \hat z = 0, 
\end{equation}
where $\hat z_j = {\rm sign} ( \hat \theta_j ) $ if $\hat \theta_j \not= 0$, $j \notin {\cal J}$,
$\hat z_j =0$ if $j \in {\cal J}$ and
$\| \hat z \|_{\infty} \le 1 $. 

We define $\theta^0$ as the solution of
$ P \dot \rho_{\theta^0} =0 $
(assumed to exist and be unique)
and assume that
$${\cal I} := I(\theta^0) := \partial P  \dot \rho_{\theta} / \partial \theta^T \vert_{\theta = \theta^0 }
$$
exists and is invertible. Its smallest eigenvalue is denoted by $\phi_0^2$. 
We define
$$ {\rm rem} (\theta - \theta^0):= {\cal I} (\theta - \theta^0) - P_n ( \dot \rho_{\theta} -
\dot \rho_{\theta^0} ). $$
The behaviour of this remainder term is studied in the next sections.

We let
$$ \Gamma_{J}:= {\cal I}_{-J, -J}^{-1} {\cal I}_{-J,J}  $$
and for all $j \in J$ and all $k \in \{ 1 , \ldots , p \}$
$$C_{J,k,j}= \begin{cases} 1 &  k =j, \cr
0  & k \not=j , \ k \in J \cr - \gamma_{J, k, j} & k \notin J \cr  \end{cases} . $$
We let $\vvert \Gamma_J \vvert_{1, {\rm off}} := \max_{j \in J} \sum_{k \in {\cal J} \backslash J } | \Gamma_{J,k,j} | $ and $\vvert C_J \vvert_{1, {\rm off}} :=1 +
\vvert \Gamma_J \vvert_{1 , {\rm off}} =
\max_{j\in J }\sum_{k \in {\cal J}  } | c_{J,k,j} | $.
We write
 $$w:= -P_n \dot \rho_{\theta^0}  . $$
 Note that if $n P_n \rho_{\theta}$ is a well-specified log-likelihood
then under regularity $n P [ w w^T ]= {\cal I}$ is the Fisher information  and $n P [C_J^T w w^T C_J ] = {\cal I}_{J,J}-
\Gamma_J^T {\cal I}_{-J, -J} \Gamma_J $. 
 
 In following theorem we derive rates for groups (Part I) and
 sharp oracle results under conditions on the worst possible
 sub-direction $\Gamma_J $ (Part II). The latter generalizes
 sharp oracle results for M-estimators presented in
\cite{van2013generic}.

\begin{theorem} \label{generalloss.theorem} We define
$$ {\cal T}_{\rm rem} := \{   \| {\rm rem} ( \hat \theta - \theta^0) \|_{\infty} \le
\lambda_{\rm rem} \} . $$
{\bf Part I}
Let
$$ \bar {\cal T}_{J, w} := \{ | \theta_J^T C_J^T w | \le\bar  \lambda_{J, \epsilon}  \sqrt {|J|} 
\sqrt {\theta_J^T ( {\cal I}_{J,J} - \Gamma_J^T {\cal I}_{-J, -J} \Gamma_J ) \theta_J } , \ 
\forall \ \theta_J \in \R^{|J|} \}.  $$
Then on $\bar {\cal T}_{J,w} \cap {\cal T}_{\rm rem}$
$$(\hat \theta_J - \theta_J^0 )^T( {\cal I}_{J,J} - \Gamma_J^T {\cal I}_{-J, -J} \Gamma_J ) 
(\hat \theta_J - \theta_J^0 ) 
  \le  \biggl (\lambda_{J,w} + \lambda \vvert C_J \vvert_{1, {\rm off}} + \lambda_{\rm rem} \vvert C_J \vvert_{1} 
  \biggr )^2
{ | J| }  /\phi_0^2 .$$
{\bf Part II} 
Let
$${\cal T}_{J,w} := \{ \| C_J^T w \|_{\infty}  \le \lambda_{J,w} \} , $$
Assume that $\vvert \Gamma_J \vvert_{1, {\rm off}} <1$ and in fact that for some
$\lambda_1 \ge 0$
$$\lambda (1- \vvert \Gamma_J \vvert_{1, {\rm off}}) >  \lambda_{J,w} + \lambda_{\rm rem}  \vvert C_J \vvert_{1 } + \lambda_1  . $$
Let $\theta \in \R^p$ be a vector with $S:= \{ \theta_j \not= 0 \} \supset {\cal J} $.
Then on ${\cal T}_{J,w} \cap {\cal T}_{\rm rem}$
$$ (\hat \theta_J - \theta_J^0 )^T( {\cal I}_{J,J} - \Gamma_J^T {\cal I}_{-J, -J} \Gamma_J ) 
(\hat \theta_J - \theta_J^0 ) $$
$$+ 
2 
\biggl [ \biggl (\lambda (1- \vvert  \Gamma_J \vvert_{1 , {\rm off}} )  -(\lambda_{ J, w}  +  \lambda_{{\rm rem}}   \vvert  C_J \vvert_1 ) - \lambda_1 \biggr ) \wedge
\lambda_1 \biggr ]  \| \hat \beta_{J }- \beta_J\|_1 $$
$$  \le 4 \lambda^2   | J \cap S |/ \phi_0^2 + ( \theta_J - \theta_J^0 )^T( {\cal I}_{J,J} - \Gamma_J^T {\cal I}_{-J, -J} \Gamma_J ) 
(\theta_J - \theta_J^0 ) .
$$

\end{theorem}

Theorem \ref{generalloss.theorem} follows using the same arguments as those used for the
proof of Theorem \ref{main.theorem}. It can be applied to obtain a global oracle inequality
($J= \{ 1 , \ldots , p \}$), a result for small coefficients ($J=S^c$) and rates for
single parameters ($J=\{ j \}$). We note however that unlike Corollary
\ref{main.corollary} all these results involve he remainder term $\lambda_{\rm rem}$
which then needs to be handled using separate arguments (e.g. applying
results from \cite{vandeG08}, \cite{pbvdg11} or \cite{Negahban12}, see also the next
section)

We also formulate an extension of Theorem
\ref{varselection.theorem}.

\begin{theorem} \label{varselection2.theorem} Let $\hat \theta$ be the unique solution of the KKT conditions
(\ref{KKT2.equation}).  Consider some set $S \supset {\cal J}$ and suppose
$\tilde \theta_S \in \Theta $ is the solution of the KKT conditions under
the restriction that the coefficients are zero outside the set $S$:
\begin{equation} \label{restricted-KKT.equation}
P_n (\dot \rho_{\tilde \theta_S})_S +   \lambda \tilde z_S =0 ,
\end{equation}
where $\tilde z_{j,S} = {\rm sign} (\tilde \theta_{j,S} )$ if $\tilde \theta_{j,S} \not= 0$, $j \in S\backslash
{\cal J} $,  $\tilde z_{j,S} = 0 $ if $j \in {\cal J}$ 
and $\| \tilde z_S \|_{\infty} \le 1 $. 
Let
${\cal T}_w := \{ \| w \|_{\infty}  \le \lambda_w\} $ and ${\cal T}_{\rm rem} :=
\{ \| {\rm rem} (\tilde \theta_S - \theta^0 ) \|_{\infty} \le \lambda_{\rm rem} \} $.
Suppose 
$$ \lambda \vvert \hat \Gamma_{S^c} \vvert_{1, {\rm off}}+ ( \lambda_{w} + \lambda_{\rm rem})
\vvert \hat \Gamma_{S^c} \vvert_1  +  \| ({\cal I}_{S^c, S^c} -
\Gamma_{S^c} {\cal I}_{S,S} 
\Gamma_{S^c} )\theta_{S^c}^0  \|_{\infty}\le   \lambda - \lambda_{w } -
\lambda_{\rm rem} 
 . $$
Then on ${\cal T}_w\cap {\cal T}_{\rm rem} $ we have $\hat \theta_j = 0 $ for all $j \notin S$. 

\end{theorem}

 \section{The remainder in terms of global norms}\label{rem.section}

\subsection{The linear model with random design}\label{rem-linear.section}

In the linear model, we write $\theta := \beta$ and we have
$$\rho_{\beta} ( X_i , Y_i ) = (Y_i - X_i \beta )^2, \ i=1 , \ldots , n , 
$$
where with some abuse of notation $X_i$ is now the $i$-th row of $X$
(i.e.\ we use the notation $X_i$ for a row and $X_j$ for a column,
the distinction only in the notation for observations and
variables: $i$ for observations and $j$ for variables).
Moreover, 
$$P_n \dot \rho_{\beta} = -X^T (Y- X \beta) /n = -X^T \epsilon  /n + \hat \Sigma (\hat \beta - \beta^0) $$
where $\hat \Sigma := X^T X / n $.
It follows that for the case of random design
$$P_n( \dot \rho_{\beta} -\dot \rho_{\beta^0}) = \hat \Sigma  (\hat \beta - \beta^0)$$
and ${\cal I} = \Sigma_0 := \EE \hat \Sigma $.
Hence
$${\rm rem} (\beta - \beta^0) = -(\hat \Sigma - \Sigma_0 )( \beta - \beta^0) $$
and
$$\| {\rm rem} (\beta - \beta^0 )\|_{\infty}  \le \| \hat \Sigma - \Sigma_0 \|_{\infty} \| \hat \beta - \beta_0 \|_1 . $$
%Taking $\lambda_{J, w} = \lambda_{J,\epsilon}$ we thus recover the result of Theorem \ref{mainrandom.theorem} with $\lambda_{J,0} := \lambda_0$ and
%${\cal T}_J \cap {\cal T}_{\rm rem}$ now defined as
%$${\cal T}_J \cap {\cal T}_{\rm rem}  := \{ \| C_J^T \epsilon \|_{\infty} /n \le \lambda_{J , \epsilon} , \
%\| \hat \Sigma - \Sigma_0 \|_{\infty} \le \lambda_0 ,\  \lambda_0 \| \hat \beta - \beta^0 \|_n \le \lambda_{\rm rem} \} .  $$
%we know that on 
%${\cal T} := \{ \| X^T \epsilon \|_{\infty} \le \lambda_{\epsilon} \}$,
%$$\| \hat \beta - \beta^0 \|_1 \le (L+1) \sqrt {s_0} \| \hat \beta - \beta^0 \|_2 , $$
%where $L= ( \lambda + \lambda_{\epsilon})/(\lambda - \lambda_{\epsilon} )$.
%It follows that for sub-Gaussian bounded design
%$${\rem} (\hat \beta - \beta^0 ) 

\subsection{Generalized linear models with random design}\label{rem-GLM.section}

We let
$$\rho_{\beta} ( X_i , Y_i) := \rho( Y_i , X_i \beta ) , \ i=1 , \ldots , n  $$
where $\{ X_i\}_{i=1}^n$ are i.i.d. $p$-dimensional row-vectors and
$\{ Y_i \}_{i=1}^n$ are i.i.d.\ response vectors with values in some
set ${\cal Y} \subseteq \R$.
Assume that for all $z, \tilde z$ 
$$| \ddot \rho_(y,z) - \ddot \rho (y, \tilde z) | \le L | z - \tilde z |  , \ \forall \ y \in {\cal Y} 
 $$
 (this can be made into a local condition)
and that $\| X \|_{\infty} \le K_X $
Then
$$\| {\rm rem} ( \beta - \beta_0)\|_{\infty}  \le K_X L \| X ( \hat \beta - \beta^0 ) \|_n^2  +
\| \hat \Sigma_{\beta^0} - \Sigma_{\beta^0} \|_{\infty} \| \hat \beta - \beta^0 \|_1. $$
Here
$$\hat \Sigma_{\beta^0} := X^T W_{\beta^0} X /n ,\  \Sigma_{\beta^0} :=
\EE \hat \Sigma_{\beta^0} $$ 
where $W_{\beta^0} := {\rm diag} ( \{ \ddot \rho ( Y_i , X_i \beta^0) \}_{i=1}^n )$. 

{\bf Standard asymptotic scenarion V} Suppose
$\ddot \rho (Y_1 , X_1 \beta^0) \ge \eta >0 $ almost surely. Assume
$K_X= {\mathcal O}(1)$, $L= {\mathcal O} (1)$ as well as
$1/ \eta = {\mathcal O} (1)$. 
Let $\lambda= {\mathcal O}_{\PP}
(\sqrt {\log p / n })$ and suppose $\| X ( \hat \beta - \beta_0 ) \|_n = 
{\mathcal O}_{\PP} ( \lambda \sqrt {s} )$ and $\| \hat \beta - \beta^0 \|_1=
o_{\PP} ( \lambda s ) $ (see e.g.\ \cite{pbvdg11}). Then
$\| {\rm rem}  (\hat \beta - \beta^0 ) \|_{\infty} = o_{\PP} (n^{-1/2})$ for
$s = o (\sqrt n / \log p ) $. 

\subsection{Exponential families}\label{rem-exponential.section}

Let $X_1 , \ldots , X_n$ be i.i.d.\ with distribution $P$.
We consider the loss function
$$ \rho_{\theta} (x) := - \sum_{j=1}^p \psi_j (x) \theta_j + d (\theta) $$
where $\theta \in \Theta$ with $\Theta$ a convex subset of $\R^p$.
Moreover, $d(\theta)$ is a twice differentiable convex function satisfying
$\dot d(\theta_0) = P\psi$. 
We assume existence of
$${\cal I} := I(\theta^0) :=\ddot d (\theta^0)   $$
and we assume that ${\cal I}^{-1}  $ exists and let $\phi_0^2$ be the smallest
eigenvalue of ${\cal I}$. 

In this case
$$w:=- P_n \dot \rho_{\theta^0} =  (P_n - P) \psi  $$
and
$${\rm rem} (\theta - \theta^0) =   {\cal I} (\theta) ( \theta - \theta^0) -[\dot d(\theta) - \dot  d(\theta_0) ].$$

%We let
%$$ \Gamma_{J}:= {\cal I}_{-J, -J}^{-1} {\cal I}_{-J,J}  $$
%and as before, for $j \in J$ and $k \in \{ 1 , \ldots , p \}$
%$$C_{J,k,j}= \begin{cases} 1 &  k =j, \cr
%0  & k \not=j , \ k \in J \cr - \gamma_{J, k, j} & k \notin J \cr  \end{cases} . $$
% We write
% $$w:= (P_n - P) \psi . $$

If we assume
\begin{equation}\label{secondderivative.equation}
\| \ddot d ( \theta ) - \ddot d (\theta^0 ) \|_{\infty} = {\mathcal O} ( \| \theta - \theta^0 \|_1 ) 
\end{equation}
then
$$\| {\rm rem} ( \theta - \theta^0) \|_{\infty} \le
\| \theta - \theta_0 \|_1 \| \ddot d ( \theta ) - \ddot d (\theta^0 ) \|_{\infty} $$
$$ = {\mathcal O} ( \|  \theta - \theta^0 \|_1^2 ) .$$

Furthermore, for exponential families
$$ d (\theta)= \log \left ( \int \exp \biggl [ \sum_{j=1}^p \psi_j \theta_j \biggr ] d \mu \right ), \ \theta \in \Theta , $$
where $\mu$ is some dominating measure for $P$. When the model
is well-specified it holds that the density $p_0:= dP / d \mu$ is equal to $p_0 = \exp[-\rho_{\theta^0} ] $.
Then (\ref{secondderivative.equation}) holds if
$\| \log p_0 \|_{\infty} = {\mathcal O} (1) $, $\max_j \| \psi_j \|_{\infty} =
{\mathcal O} (1) $ and $\| \theta - \theta^0 \|_1 = {\mathcal O} (1)$.
 
{\bf Standard asymptotic scenario VI} Suppose that $\lambda = {\mathcal O}_{\PP} (\sqrt {\log p / n })$ and that
$\| \hat \theta - \theta^0 \|_1 = {\mathcal O}_{\PP} ( \lambda s )$ and 
$\| \hat \theta - \theta^0 \|_2 = {\mathcal O}_{\PP} ( \lambda \sqrt s )$
(see e.g.\ \cite{van2013generic} for such results for general high-dimensional
models). Then 
$\| {\rm rem } (\hat \theta - \theta^0 ) \|_{\infty} = o_{\PP} (n^{-1/2})$ for
$s = o( n^{1/4} / \log p ) $.

\section{Brouwer's fixed point theorem}\label{Brouwer.section}

The remainder term ${\rm rem} (\hat \theta - \theta^0)$ will generally only be small if
$\hat \theta$ is close enough to $\theta^0$. 
If the global rate of convergence is too slow we need a technique
different from the one of the previous section to deal with the remainder term.
Here, Brouwer's fixed point theorem can be useful. The idea is from \cite{ravikumar2011high}.

\begin{lemma} \label{Brouwer.lemma} Let $\hat \theta$ be the unique solution
in $\R^p$ 
of the estimating equations
$$   \hat \theta - \theta^0 = G( {\rm rem } (\hat \theta - \theta^0)) + v (\hat \theta - \theta^0) +u_0 $$
where $G  :\R^p \rightarrow \R^p$ and 
$v :\ \R^p \rightarrow \R^p$ are functions and $u_0 \in \R^p$ is
a constant vector.  Let $q \ge 1$.
Suppose that for some constants $K$ and $\varepsilon$
$$G (B_q (1) ) \in B_q (K ) ,$$
$$\|  u_0 \|_q \le K \varepsilon, \ 
\sup_{\theta} \| v (\theta - \theta^0) \|_q \le K \varepsilon $$
and that
$$  \| {\rm rem} ( B_q (3K \epsilon) ) \|_q   \le \varepsilon . $$
Then $\| \hat \theta - \theta^0 \|_q \le 3K \varepsilon$.
\end{lemma}

\section{The irrepresentable condition and rates in sup-norm}\label{irrepresentable.section}
We showed in Theorems
\ref{varselection.theorem} and \ref{varselection2.theorem} that an
irrepresentable condition on the set of small variables can be used to show
that these parameters are estimated as being zero.
However, 
in some cases there are no ready-to-use results
to handle the remainder term. 
One may then need rather strong conditions to deal with this. 
We now show that irrepresentable conditions can lead to
convergence in sup-norm. The remainder term
then needs to be small for sup-norm neighbourhoods which are
smaller than the $\ell_1$- or $\ell_2$-neighbourhoods these imply.
The idea is as in \cite{ravikumar2011high} but with the extension that
we only need approximate sparsity (i.e. there may be many
non-zero but small coefficients).

In the next theorem we use the notation of Section
\ref{general.section},  and assume that $\hat \theta$ is the unique solution of the KKT conditions
(\ref{KKT2.equation}). The proof is based on Lemma
\ref{Brouwer.lemma}.

\begin{theorem} \label{varselection3.theorem}
%Let $\hat \theta $ be the unique solution of 
%$$ P_n \dot \rho_{\hat \theta} +   \lambda \hat z =0 ,$$
%where $\hat z_{j} = {\rm sign} (\hat \theta_{j} )$ if $\hat \theta_{j} \not= 0$, $j = 1 , \ldots, p $, 
%and $\| \hat z \|_{\infty} \le 1 $.
Let $S \supset {\cal J}$ and suppose the solution $\tilde \theta_S \in \Theta$ of the restricted KKT
conditions (\ref{restricted-KKT.equation}) exists.
Assume for some $\kappa_{S^c} \ge 1$
$$\vvert {\cal I}_{S^c, S^c}- \Gamma_{S^c}^T {\cal I}_{S,S} \Gamma_{S^c}   \vvert_1 \le
\kappa_{S^c}, $$
for some $K_S \ge 1$
$$\vvert {\cal I}_{S,S}^{-1} \vvert_1 \le K_S$$
and that for some $\lambda_{\rm thres}$
$$ \| \theta_{S^c}^0 \|_{\infty} \le \lambda_{\rm thres} .$$
Let
$${\cal T}_w := \{ \| w \|_{\infty} \le \lambda_w  \} , $$
$${\cal T}_{\rm rem} := \{ \| {\rm rem} ( \theta_S - \theta^0 )\|_{\infty} \le
\lambda_{\rm rem}  \ \forall \ \theta_S \in \Theta:  \  \| \theta_S - \theta^0 \|_{\infty} \le 3K ( \lambda_w + \lambda ) \} . $$
Suppose 
$$
\lambda \vvert \Gamma_{S^c} \vvert_{1, {\rm off}}+ (\lambda_w + \lambda_{\rm rem})  \vvert \Gamma_{S^c} \vvert_{1} +
\kappa_{S^c} \lambda_{\rm thres} \le \lambda - \lambda_w - \lambda_{\rm rem}   . $$
Then on ${\cal T}_w \cap {\cal T}_{\rm rem} $, we have
$\| \hat \theta_S - \theta_S^0 \|_{\infty} \le 3K_S ( \lambda_w + \lambda)$ and
$\| \hat \theta_{S^c} \|_{\infty} =0 $. 
\end{theorem}

\section{Estimating a precision matrix}\label{precision.section}

We investigate the remainder term for the graphical Lasso.
The approach is again similar to \cite{ravikumar2011high}. Our main extension is that
 we no longer assume that the truth is exactly sparse (in the sense of having many
 parameters exactly equal to zero) but only approximately sparse
 (i.e. having a sparse approximation).  This extension may be important when
 applying the results for obtaining confidence intervals (see \cite{Jankova2014})
 because approximate sparsity appears more in line with the concept of honesty for confidence
 intervals.
 
 Let $X$ be an $n \times p$ matrix with i.i.d.\ rows with distribution $P$.
 We let $\hat \Sigma := X^T X / n $ and $\Sigma_0 := \EE \hat \Sigma$.
 We assume that $\Theta_0 := \Sigma_0^{-1} $ exists. The matrix
 $\Theta_0$ is called the precision matrix.
 We consider the estimator
 $$\hat \Theta = \arg \min_{\Theta \ {\rm p.s.d}}
 {\rm trace} ( \hat \Sigma \Theta )- \log {\rm det} (\Theta) + 2 \lambda \| \Theta \|_{1, {\rm off}} , $$
 where $\| \Theta \|_{1 , {\rm off} } := \sum_j \sum_{k \not= j} | \Theta_{k,j} | $. 
 (We note that $\Theta$ is now the parameter (not the parameter space)
 and the parameter space is the set of all positive semi-definite (p.s.d.) matrices.)
 Observe this corresponds to a loss function from the exponential family, namely 
 $$\rho_{\Theta} (x) =- \sum_{j,k}\psi_{k,j} (x)\Theta_{k,j}  + d (\Theta), \ \psi_{k,j} (x)=
 - x_k x_j , \ d (\Theta) =-
 \log {\rm det} ( \Theta) . $$

Consider the KKT-conditions
$$\hat \Sigma - \hat \Theta^{-1} + \lambda \hat Z =0 $$
where $\hat Z_{j,j} =0$ for all $j$ and for $j \not= k $, $\hat Z_{j,k} := {\rm sign} (\hat \Theta_{j,k} )$
when $\hat \Theta_{j,k} \not= 0 $. Moreover $\| \hat Z \|_{\infty} \le 1 $. 
We define $W = \hat \Sigma - \Sigma_0 $ and  write this as
$$W - ( \hat \Theta^{-1} - \Theta_0^{-1} ) + \lambda \hat Z = 0 . $$

Note that for $\hat \Delta:= \hat \Theta - \Theta_0$
$$-( \hat \Theta^{-1} - \Theta_0^{-1} ) = \hat \Theta^{-1} \hat  \Delta \Theta_0^{-1} $$
$$:= \Theta_0^{-1} \hat \Delta \Theta_0^{-1} - {\rm rem} (\hat \Delta), $$
where
$$ {\rm rem}  (\Delta) := \biggl ( (\Theta_0 + \Delta)^{-1} - \Theta_0^{-1}+ \Theta_0^{-1} \Delta \Theta_0^{-1}\biggr ) $$
is the remainder term.

We let $\Theta_S$ be a symmetric positive definite matrix which has zeroes outside the set $S$. 
Moreover, $\Theta_S^0$ denotes the matrix $\Theta_0$ with all entries in $S^c$
set to zero (this matrix may not be positive semi-definite) and $\Theta_{S^c}^0 := \Theta_0 -
\Theta_S^0 $.
We moreover write
$$ \vvert \Theta_S^0  \vvert_0 := \max_{j} | \{ k: \  (k,j) \in S ,\ \Theta_{0, k,j} \not= 0 \}|  . $$

\begin{lemma}\label{graphicalrem.lemma}
Suppose that for some $\epsilon_0$ and $\eta_0$
$$ 
\vvert \Theta_0^{-1} \vvert_{1} \biggl ( \varepsilon_0   \vvert \Theta_S^0 \vvert_0 + \vvert \Theta_{ S^c}^0 \vvert_1 \biggr ) \le \eta_0 <1 . $$
Then for $\| \Theta_S - \Theta_0 \|_{\infty} \le \epsilon_0 $ we have
$${\rm rem} ( \Theta_S - \Theta_0) \le \eta_0 \epsilon_0 \vvert \Theta_0^{-1} \vvert_{1}^2/(1- \eta_0) . $$

\end{lemma}

{\bf Standard asymptotic scenario VII} We combine Lemma \ref{graphicalrem.lemma}
with Theorem \ref{varselection3.theorem}, where in the latter we use the vectorized versions
of the matrix parameter $\Theta$ to define the matrices involved. 
 Suppose that $\vvert \Theta_0^{-1} \vvert_1 =
{\mathcal O} (1)$ (recall $\Theta_0^{-1} = \Sigma_0 $), $\vvert {\cal I}_{S^c, S^c} - \Gamma_{S^c}^T {\cal I}_{S,S} \Gamma_{S^c} \vvert_1=
{\mathcal O} (1)$ and $\vvert {\cal I}_{S,S}^{-1} \vvert_1 = {\mathcal O}(1)$. 
Assume that the probability of ${\cal T}_w$ (with $w$ the vectorized version of $W$) tends to one
for $\lambda_w= {\mathcal O} ( \sqrt {\log p / n } )$
 and let $\lambda = {\mathcal O}_{\PP} ( \sqrt {\log p / n} )$ be suitably large.
 Let $S := \{ j,k:\ | \Theta_{0,j,k} | > \lambda_{\rm thres} \}$ where
 $\lambda_{\rm thres}= {\mathcal O} ( \sqrt {\log p / n })$ is suitably small.
Assume moreover $\vvert \Theta_S^0 \vvert_1 = d$ (this is up to
small coefficients the maximal edge degree of the matrix $\Theta_0$) 
with $d = O ( \sqrt { n } / \log p )$ and that
$\vvert \Theta_{S^c}^0 \vvert_1 = {\mathcal O} ( \lambda d ) $.
Finally assume that $\vvert \Gamma_{S^c} \vvert_1 < 1$ is sufficiently small (but not
necessarily tending to zero). Then $\| \hat \Theta - \Theta_0 \|_{\infty} =
{\mathcal O}_{\PP} ( \sqrt {\log p / n } )$ and moreover
$\hat \Theta_{j,k} = 0 $ for all $(j,k) \in S $. 
An example of an approximately sparse case where this result can be applied is
where
the columns of $\Theta_0$ have a uniformly finite $\ell_r$-``norm" for some $0 < r < 1$, i.e.\ where
$$\max_{j} \sum_{k } | \Theta_{0,k,j} |^r = {\mathcal O} (1) .$$

\section{Proofs} \label{proofs.section}

\subsection{Proofs for Section \ref{fixed-design.section}: The linear model with fixed design}

We first proof the results for a single parameter.

{\bf Proof of Lemma \ref{oneparameter.lemma}.} 
Define for $t \in \R$, $\hat \beta (t) := \hat \beta + t \hat c_j  $. By the KKT conditions
(see (\ref{KKT.equation}))
the sub-gradient of 
$\| Y - X \hat \beta (t) \|_n^2/2 +  \lambda \| \beta(t) \|_1 $ at $t =0$ is equal to zero, i.e.
$$- (\hat c_j  X )^T (Y- X \hat \beta ) /n +  \lambda   \hat c_j^T \hat z  =0  $$
where $\hat z_j = {\rm sign} (\hat \beta_j )$ if $\hat \beta_j \not= 0$ and
moreover $\| \hat z \|_{\infty} \le 1$. 
This can be rewritten to
$$-(X_j - X_{-j} \hat \gamma_j)^T \epsilon/n +
(X_j - X_{-j} \hat \gamma_j)^T X_j ( \hat \beta_j - \beta_j^0) /n $$ $$ +
(X_j - X_{-j} \hat \gamma_j)^T X_{-j} ( \hat \beta_{-j} - \beta_{-j}^0)/n
=\lambda \hat c_j^T \hat  z  . $$
But from the KKT conditions, using the counterpart of (\ref{useful.equation}) for $\hat \gamma_j$,
$$ (X_j - X_{-j} \hat \gamma_j)^T X_j /n= \tilde  \tau_j^2  .$$
Moreover, again by the KKT conditions, using the counterpart of 
(\ref{KKT.equation})  for $\hat \gamma_j$,
$$ \| (X_j - X_{-j} \hat \gamma_j)^T X_{-j}\|_{\infty}/n  \le \lambda_j .$$
It is also clear that
%$$ | (X_j - X_{-j} \hat \gamma_j)^T \epsilon |/n \le \lambda_{\epsilon}  \hat \tau_j  $$ and
$$|\hat c_j^T\hat z  | \le  \| \hat c_j \|_1 . $$
We thus find on ${\cal T}_{ j , \epsilon}  $ that
$$\tilde \tau_j^2 | \hat \beta_j - \beta_j^0 | \le  \lambda_{j, \epsilon} \hat \tau_j +
\lambda_j \| \hat \beta_{-j} - \beta_{-j}^0 \|_1 + \lambda  \| \hat c_j \|_1  . $$
Whence the first result.
The second result follows from
$$( \hat \beta_j - \beta_j^0 ) \hat c_j^T \hat z \ge |\hat \beta_j | - |\beta_j^0| -
\| \hat \gamma_j \|_1 | \hat \beta_j - \beta_j^0 | . $$
 \hfill $\sqcup \mkern -12mu \sqcap$
 
 {\bf Proof of Lemma \ref{all.lemma}.} This follows immediately from Corollary \ref{oneparameter.corollary}. 
\hfill $\sqcup \mkern -12mu \sqcap$

{\bf Proof of Lemma \ref{de-sparsify.lemma}.} Let $\hat \Sigma := X^T X / n $ be the Gram matrix.
By the KKT conditions for $\{ \hat \gamma_j \}$, we have
$$ \hat \Sigma  \hat \Theta -I = \lambda_0 \hat Z \hat T \tilde  T^{-2}  , $$
where $\hat Z_{j,j} = 0 $ for all $j$ and for $k \not= j$,  $\hat z_{k,j} = {\rm sign} (\hat \theta_{k,j} )$
when $\hat \theta_{k,j}
\not= 0$. Moreover $\| \hat Z \|_{\infty} \le 1$.
Therefore
$$\hat T^{-1} \tilde T^2  ( \hat b - \beta^0) = \hat T^{-1} \tilde T^2  (\hat \beta - \beta^0)  +
\hat T^{-1}  \tilde T^2  \hat \Theta^T X^T ( Y- X \hat \beta )  /n$$
 $$ = \hat T^{-1}  \tilde T^2 \hat \Theta^T X^T \epsilon/n - \hat T^{-1} \tilde T^2  (\hat \Theta^T \hat \Sigma -I) (\hat \beta - \beta^0)  $$ $$=  \hat T^{-1} \tilde T^2  \hat \Theta^T X^T \epsilon/n -
 \lambda_0 \hat Z^T ( \hat \beta - \beta^0) . $$
 Clearly, $\| \hat Z^T ( \hat \beta - \beta^0) \|_{\infty} \le \| \hat Z \|_{\infty} 
\|  \hat \beta - \beta^0 \|_1  \le \| \hat \beta - \beta^0 \|_1 $.
\hfill $\sqcup \mkern -12mu \sqcap$

{\bf Proof of lemma \ref{Dantzig.lemma}.} 
 As in the proof of Lemma \ref{oneparameter.lemma}
$$(X_j - X_{-j} \hat \gamma_j)^T \epsilon/n +
(X_j - X_{-j} \hat \gamma_j)^T X_j ( \hat \beta_j - \beta_j^0) /n $$ $$ +
(X_j - X_{-j} \hat \gamma_j)^T X_{-j} ( \hat \beta_{-j} - \beta_{-j}^0)/n
=\lambda \hat  z^T \hat c_j . $$
Now by definition
$$ | (X_j - X_{-j} \hat \gamma_j)^T X_{-j} ( \hat \beta_{-j} - \beta_{-j}^0)|/n \le
\lambda_j \| \hat \beta_{-j} - \beta_{-j}^0 \|_1 . $$
Furthermore
$$ (X_j - X_{-j} \hat \gamma_j)^T X_j  /n =\hat \tau_j^2 + (X_j - X_{-j} \hat \gamma_j)^T X_{-j} \hat \gamma_j /n
\ge \hat \tau_j^2 - \lambda_j \| \hat \gamma_j \|_1 . $$
Continuing as in the proof of Lemma \ref{oneparameter.lemma}, we get on ${\cal T}_{j, \epsilon} $
$$\hat \tau_j^2 | \hat \beta_j - \beta_j^0 | \le \lambda_{j, \epsilon} \hat \tau_j +
\lambda_j \| \hat \beta_{-j} - \beta_{-j}^0 \|_1 + \lambda  \| \hat c_j \|_1 +  
\lambda_j \| \hat \gamma_j \|_1 | \hat \beta_j - \beta_j^0 |  .$$
\hfill $\sqcup \mkern -12mu \sqcap$

The main result for a group of parameters has a somewhat more involved proof.

{\bf Proof of Theorem \ref{main.theorem}.} 
Recall the KKT conditions
$$
-X^T (Y- X \hat \beta )/n + \lambda \hat z =0 $$
where $\hat z_j = {\rm sign} (\hat \beta_j ) $ if $\hat \beta_j \not=0 $ and
$\| \hat z \|_{\infty} \le 1 $. 
Multiply by $( \hat \beta_J - \beta_J)^T \hat C_J^T$
to find that
$$- ( \hat \beta_J - \beta_J)^T \hat C_J^T X^T ( Y- X \hat \beta ) / n +
\lambda ( \hat \beta_J- \beta_J)^T \hat C_J^T \hat z = 0 . $$
Rewrite this to
$$-(\hat \beta_J - \beta_J)^T \hat C_J^T X^T \epsilon /n + 
(\hat \beta_J - \beta_J)^T \hat C_J^T X^T X \hat C_J (\hat \beta_J - \beta_J^0)  /n$$ $$+
(\hat \beta_J - \beta_J)^T\hat C_J^T X^T X_{-J} \biggl (( \hat \beta_{-J} - \beta_{-J}^0 )+ \hat \Gamma_J ( \hat \beta_J - \beta_J^0)\biggr )/n $$ $$ + \lambda \| \hat \beta_J \|_1 - \lambda (\beta_J)^T \hat z_J 
- \lambda  (\hat \beta_J - \beta_J)^T \hat \Gamma_J^T\hat z_{-J}  =0. $$

By the KKT conditions for $\hat \Gamma_J$
$$\|  (X \hat C_J)^T X_{-J}  \|_{\infty}/n \le \lambda_0 . $$

We therefore get
\begin{equation}\label{continue.equation}
(\hat \beta_J - \beta_J)^T \hat C_J^T X^T X \hat C_J ( \hat \beta_J - \beta_J^0 )  /n 
\le  
| ( \hat \beta_J - \beta_J)^T \hat C_J^T X \epsilon | /n 
\end{equation} 
$$+\biggl (\lambda_{J,0} \| \hat \beta_{-J}- \beta_{-J}^0 \|_1 +  \lambda \vvert \hat \Gamma_J \vvert_1  \biggr )  \| \hat \beta_J - \beta_J \|_1
   $$
$$  +  \lambda_{J,0}  \vvert \hat \Gamma_J \vvert_1 \| \hat \beta_J - \beta_J \|_1^2 +
 \lambda \| \beta_J \|_1   
  -  \lambda 
\| \hat \beta_J \|_1 .$$
{\bf Part I} 
To obtain the first equation, we apply (\ref{continue.equation}) with $\beta_J = \beta_J^0$ and we
use that $$\| \hat \beta_J - \beta_J^0 \|_1^2 \le | J| \| X \hat C_J ( \hat \beta_J - \beta_J^0 ) \|_n^2/
\hat \phi_J^2  . $$
Then we get  on $\bar {\cal T}_{J, \epsilon}$
$$\biggl ( 1- { \lambda_{J,0} |J| \vvert \hat \Gamma_J \vvert_1 \over
\hat \phi_J^2 }  \biggr ) \| X \hat C_J ( \hat \beta_J - \beta_J^0 ) \|_n^2
 \le   \bar \lambda_{ J,\epsilon} \sqrt {|J|} \| X \hat C_J ( \hat \beta_J -
\beta_J^0 ) \|_n $$ $$+  \biggl ( \lambda_{J,0} \| \hat \beta_{-J}-
\beta_{-J}^0 \|_1 
+ \lambda \vvert \hat \Gamma_J \vvert_1  \biggr )  \| \hat \beta_J - \beta_J^0 \|_1 + 
 \lambda \| \beta_J^0 \|_1   
  -  \lambda 
\| \hat \beta_J \|_1 $$
$$ \le   \bar \lambda_{ J,\epsilon}  \sqrt {|J|} \| X \hat C_J ( \hat \beta_J -
\beta_J^0 ) \|_n + \biggl ( \lambda_{J,0} \| \hat \beta_{-J}-
\beta_{-J}^0 \|_1 
+ \lambda \vvert \hat C_J \vvert_1  \biggr )  \| \hat \beta_J - \beta_J^0 \|_1
$$
$$ \le {\sqrt{ |J| } \over \hat \phi_J} \biggl (\bar  \lambda_{ J,\epsilon} \hat \phi_J  +  \lambda_{J,0} \| \hat \beta_{-J}-
\beta_{-J}^0 \|_1 
+ \lambda \vvert \hat C_J \vvert_1  \biggr ) \| X \hat C_J ( \hat \beta_J - \beta_J^0 ) \|_n.$$
This gives the first result (\ref{for-desparsifying.equation}).

The second result (\ref{simpler.equation}) follows from (\ref{continue.equation}) by similar arguments
and the bounds $\| \hat \beta_J - \beta_J^0 \|_1 \le 
\| \hat \beta - \beta^0 \|_1$ and $\| \hat \beta_{-J}- \beta_{-J}^0 \|_1 \le \| \hat \beta - \beta^0 \|_1$. 
%  $$ \le  { \sqrt {J} \over \hat \phi( \vvert \hat \Gamma_J \vvert_1 , J) } \biggl (  \lambda_{ \epsilon}+
%   \lambda+ 
%  \lambda_0 \| \hat \beta - \beta^0 \|_1   \biggr )  \vvert \hat C_J \vvert_1
% \| X \hat C_J ( \hat \beta_J - \beta_J^0 ) \|_n  . $$
% This gives the first result.

{\bf Part II} 
We may also use (\ref{continue.equation}) to obtain that on ${\cal T}_{J, \epsilon}$
\begin{equation} \label{restart.equation}
(\hat \beta_J - \beta_J)^T \hat C_J^T X^T X  \hat C_J ( \hat \beta_J - \beta_J^0 )   /n 
\end{equation}
$$ \le    \biggl [ \lambda_{J, \epsilon} 
 +  (\lambda_{J,0} \| \hat \beta-
\beta^0 \|_1+ \lambda ) \vvert \hat C_J \vvert_1 \biggr ] 
\|\hat \beta_{J \cap S}  - \beta_{J \cap S}  \|_1 $$ $$+  \biggl [ \lambda_{ J, \epsilon}  +  \lambda_{J,0} \| \hat \beta -
\beta^0 \|_1  \vvert \hat C_J \vvert_1+ \lambda \vvert \hat \Gamma_J  \vvert_1 \biggr ] \| \hat \beta_{J \backslash S} \|_1  
  -  \lambda \| \hat \beta_{J \backslash S}  \|_1. $$
  
  If 
  $$ \| X \hat C_J ( \hat \beta_J - \beta_J^0 ) \|_n^2 + 2 \lambda_1 \| \hat \beta_J - \beta_J \|_1  \le
  \| X \hat C_J ( \beta_J - \beta_J^0 ) \|_n^2 $$
  we are done. So let us assume in the rest of the proof that
   $$ \| X \hat C_J ( \hat \beta_J - \beta_J^0 ) \|_n^2 + 2 \lambda_1 \| \hat \beta_J - \beta_J \|_1  \ge
  \| X \hat C_J ( \beta_J - \beta_J^0 ) \|_n^2 .$$
  Then we have
  $$(\hat \beta_J - \beta_J)^T \hat C_J^T X^T X \hat C_J (\hat \beta_J - \beta_J^0) /n 
    = {1 \over 2} \| X\hat C_J (\hat \beta_J - \beta_J^0 ) \|_n^2 $$
  $$+{1 \over 2}
  \| X\hat C_J (\hat \beta_J - \beta_J ) \|_n^2 - {1 \over 2} \| X\hat C_J (\beta_J - \beta_J^0 ) \|_n^2$$
 $$ \ge  {1 \over 2} \| X\hat C_J (\hat \beta_J - \beta_J ) \|_n^2 - \lambda_1 \| \hat \beta_J - \beta_J \|_1 
 \ge - \lambda_1 \| \hat \beta_J - \beta_J \|_1 .$$
  This gives by (\ref{restart.equation})
$$
\biggl [ \lambda -  \lambda_{J, \epsilon} -  \lambda_{J,0} \| \hat \beta-
\beta^0 \|_1  \vvert \hat C_J \vvert_1- \lambda \vvert \hat \Gamma_J  \vvert_1 - \lambda_1 \biggr ] \| \hat \beta_{J \backslash S}  \|_1
$$
 $$\le \biggl  [ \lambda_{J, \epsilon} +  (\lambda_{J,0} \| \hat \beta-
\beta^0 \|_1+ \lambda)\vvert \hat C_J \vvert_1 + \lambda_1 \biggr ] 
\|\hat \beta_{J \cap S}  - \beta_{J \cap S} \|_1  , $$
and hence on ${\cal T}_{J, \epsilon} \cap {\cal T}_{J, {\rm rem}}$
$$ \biggl ( \lambda (1- \vvert \hat \Gamma_J \vvert_1 )-
( \lambda_{J, \epsilon} + \lambda_{J,{\rm rem}}   \vvert \hat C_J \vvert_1 ) - \lambda_1  \biggr ) \|  \hat \beta_{J \backslash S} \|_1   $$ $$ \le
(\lambda_{J, \epsilon} + (\lambda + \lambda_{J,{\rm rem} }  ) \vvert \hat C_J \vvert_1 + \lambda_1 )  \|\hat \beta_{J \cap S}  - \beta_{J \cap S}\|_1 .$$
Moreover
$$ \| \hat \Gamma_J ( \hat \beta_J - \beta_J ) \|_1 \le
\vvert \hat \Gamma_J \vvert_1 \| \hat \beta_J - \beta_J \|_1 =
\vvert \hat \Gamma_J \vvert_1 \| \hat \beta_{J \backslash S}  \|_1 + 
\vvert \hat \Gamma_J \vvert_1 \|\hat \beta_{J \cap S}  - \beta_{J \cap S}  \|_1  .$$
Thus
$$\| \hat \beta_{J \backslash S}  \|_1 + \| \hat \Gamma_J ( \hat \beta_J - \beta_J ) \|_1  \le
(1+ \vvert \hat \Gamma_J \vvert_1)  \| \hat \beta_{J \backslash S}  \|_1 + 
\vvert \hat \Gamma_J \vvert_1 \|\hat \beta_{J \cap S}  - \beta_{J \cap S} \|_1 $$
$$ \le L  \|\hat \beta_{J \cap S}  - \beta_{J \cap S}  \|_1 . $$
It follows that
$$\| \hat \beta_{J \cap S} - \beta_{J \cap S}\|_1 \le
{\sqrt { | J \cap S | } \over \hat \phi (L, J \cap S) } \| X \hat C_J ( \hat \beta_J - \beta_J ) \|_n. $$
But then
\begin{equation}\label{3par.equation}
{1 \over 2} \| X\hat C_J (\hat \beta_J - \beta_J^0 ) \|_n^2 
  +{1 \over 2}
  \| X\hat C_J (\hat \beta_J - \beta_J ) \|_n^2 - {1 \over 2} \| X\hat C_J (\beta_J - \beta_J^0 ) \|_n^2
\end{equation}
 $$+  \biggl [ \lambda (1- \vvert \hat \Gamma_J \vvert_1 ) - (\lambda_{ J, \epsilon}  +  \lambda_{J,{\rm rem}}   \vvert \hat C_J \vvert_1 ) - \lambda_1 \biggr ] \| \hat \beta_J - \beta_J \|_1  
 $$ $$
 \le    2\lambda 
{\sqrt { | J \cap S | } \over \hat \phi (L, J \cap S) } \| X \hat C_J ( \hat \beta_J - \beta_J ) \|_n
$$ $$ \le 2 \lambda^2 
 {{ | J \cap S | } \over  \hat \phi^2  (L, J \cap S) } +  {1 \over 2} \| X \hat C_J ( \hat \beta_J - \beta_J ) \|_n^2 .$$
This gives
$$\| X \hat C_J ( \hat \beta_J - \beta_J^0 ) \|_n^2 +2 \biggl [ \lambda (1- \vvert \hat \Gamma_J \vvert_1 )  -(\lambda_{ J, \epsilon}  +  \lambda_{J,{\rm rem}}   \vvert \hat C_J \vvert_1 )- \lambda_1 \biggr ] \| \hat \beta_{J }- \beta_J\ |_1  $$ $$
  \le {{4 \lambda^2  | J \cap S | } \over \hat \phi ^2 (L, J \cap S) }
 + \| X \hat C_J ( \beta_J - \beta_J^0) \|_n^2  . $$

\hfill $\sqcup \mkern -12mu \sqcap$

We now show that under certain conditions the estimator puts values in $S^c$ to zero.

{\bf Proof of Theorem \ref{varselection.theorem}.}
Let $\tilde \beta_S$ be the solution of the KKT conditions under
the restriction that the coefficients are zero outside the set $S$:
$$-X_S^T (Y - X_S \tilde \beta_S)/n  + \lambda \tilde z_S =0, $$
where for $j \in S$, $\tilde z_{j, S}= {\rm sign} ( \tilde \beta_{j,S} )$
if $\tilde \beta_{j,S} \not= 0 $ and where $\|\tilde z_S \|_{\infty} \le 1 $.
Define
$$\lambda \tilde z_{S^c} := X_{S^c}^T ( Y- X_S \tilde \beta_S)/n  .
$$
Then
$$\lambda \tilde z_{S^c}
=
\lambda \hat \Gamma_{S^c}^T \tilde z_S +
\hat C_{S^c}^T X^T X \hat C_{S^c} \beta_{S^c}^0/n  +\hat C_{S^c} ^T X^T \epsilon/n  . $$
It follows that on ${\cal T}_{\epsilon}$
$$ \| \tilde  z_{S^c} \|_{\infty} \le 1 . $$
Hence $\tilde \beta_S$ is a solution of the unrestricted KKT
conditions (\ref{KKT.equation}) and hence $\tilde \beta_S = \hat \beta$. 

\hfill $\sqcup \mkern -12mu \sqcap$

The final proof of this section concerns the de-sparsification of a group of variables.

{\bf Proof of Lemma \ref{desparsify-group.lemma}.} We have
%$$\hat C_J^T X^T X_J \hat b_J/n =
%\hat C_J^T X^T X_J  \hat \beta_J /n+ \hat C_J^T X^T (Y - X \hat \beta )/n $$
%$$= \hat C_J^T X^T X_J  \hat \beta_J /n + \hat C_J X^T \epsilon/n
%- \hat C^T X^T X_J ( \hat \beta_J - \beta_J^0) -
%\hat C^T X^T X_{-J} ( \hat \beta_{-J} - \beta_{-j}^0) $$
%Hence
$$\hat C_J^T X^T X_J (\hat b_J- \beta_J^0)/n=
\hat C_J^T X^T X_J ( \hat \beta_J- \beta_J^0) /n $$ $$+ \hat C_J X^T \epsilon/n
- \hat C^T X^T X_J ( \hat \beta_J - \beta_J^0)/n -
\hat C^T X^T X_{-J} ( \hat \beta_{-J} - \beta_{-J}^0) /n$$
$$ = \hat C_J X^T \epsilon/n - \hat C^T X^T X_{-J} ( \hat \beta_{-J} - \beta_{-J}^0) /n.$$
The result now follows using the same arguments as for deriving
(\ref{for-desparsifying.equation}) and (\ref{simpler.equation}) of Theorem
\ref{main.theorem}.

%But then
%$$\| X \hat C_J ( \hat b_J - \beta_J^0) \|_n^2 =
% + ( \hat b_J - \beta_J^0)^T \hat C_J X^T \epsilon/n + 
% ( \hat b_J - \beta_J^0)^T \hat C^T X^T X_{-J} ( \hat \beta_{-J} - \beta_{-J}^0) /n$$
%$$ + ( \hat b_J - \beta_J^0)^T \hat C_J^T X^T X_{-J} \hat \Gamma_J ( \hat \beta_J  - \beta_J^0)/n . $$
%$$ \le \sqrt {|J|} \| X \hat C_J ( \hat b_J - \beta_J^0) \|_n \lambda_{J,\epsilon} +
%\lambda_{J,0} \| \hat \beta_{-J}- \beta_{-J}^0 \|_1  \| \hat b_J - \beta_J^0 \|_1 $$ $$
%+ \lambda_{J,0} \vvert \hat \Gamma_J \vvert_1 \| \hat b_J - \beta_J^0 \|_1
%\| \hat \beta_J - \beta_J^0 \|_1  $$
%$$ \le \sqrt {|J|} \| X \hat C_J ( \hat b_J - \beta_J^0) \|_n \lambda_{J,\epsilon} +
%(\lambda_{J,0}  \| \hat \beta- \beta^0 \|_1 ) \| X \hat C_J ( \hat b_J -
%\beta_J^0) \|_n  \sqrt {|J|}/ \hat \phi_J $$ $$+ ( \lambda_{J, \epsilon} +
%( \lambda + \lambda_{J,0} \| \hat \beta- \beta^0 \|_1 ) \vvert \hat C_J \vvert_1 )
% \| X \hat C_J ( \hat b_J -
%\beta_J^0) \|_n  \lambda_{J,0}\vvert \hat \Gamma_J \vvert_1 { |J| } / \hat \phi_J^2  $$
%where we used the second bound of Theorem \ref{main.theorem} for
%$\| \hat \beta_J - \beta_J^0 \|_1$

\hfill $\sqcup \mkern -12mu \sqcap$

\subsection{Proof for Section \ref{random.section}: Random design.}

{\bf Proof of Theorem \ref{mainrandom.theorem}.}
This follows using the same arguments as in the proof of Theorem \ref{main.theorem}.
\hfill $\sqcup \mkern -12mu \sqcap$

\subsection{Proof for Section \ref{general.section}: General loss functions}

{\bf Proof of Theorem \ref{generalloss.theorem}.}
By the KKT conditions (\ref{KKT2.equation})
$$-w + {\cal I} ( \hat \theta - \theta^0) -   {\rm rem}  (\hat \theta - \theta^0) + \lambda \hat z=0 . $$
Multiplying by $(\hat \theta_J - \theta_{J} )^T C_J^T$ and rewriting gives
$$-(\hat \theta_J - \theta_J)^T C_J^Tw
 + 
(\hat \theta_J - \theta_J )({\cal I}_{J,J} - \Gamma_J^T {\cal I}_{-J, -J } \Gamma_J) (\hat \theta_J - \theta_J^0) $$
$$ + \lambda (\hat \theta_J - \theta_J)^T \hat z_{J} +
\lambda ( \hat \theta_J - \theta_J^0)^T \Gamma_J^T \hat z_{-J} $$
$$ -  (\hat \theta_J - \theta_J)^T
  {\rm rem}_{-J}  (\hat \theta-\theta^0 ) - ( \hat \theta_J - \theta_J^0)^T \Gamma_J^T    {\rm rem}_{-J}  (\hat \theta-\theta^0) = 0 . $$
  The results now follow in the same manner as for Theorem \ref{main.theorem}.

\hfill $\sqcup \mkern -12mu \sqcap$

{\bf Proof of Theorem \ref{varselection2.theorem}.} This follows from
the same arguments as those used for proving Theorem
\ref{varselection.theorem}. \hfill $\sqcup \mkern -12mu \sqcap$

\subsection{Proof for Section \ref{Brouwer.section}: Brouwer's fixed point theorem}

{\bf Proof of Lemma \ref{Brouwer.lemma}.} 
Let $F(\delta) := G ( {\rm rem } (\delta) ) + v (\delta) +u_0$.
Then for $\| \delta \|_q \le  3K \varepsilon $
$$\| F(\delta ) \|_q \le \| G ({\rm rem} (\delta )) \|_q + \| v (\delta ) \|_q + \|
 u_0 \|_q $$
$$ \le K \| {\rm rem} (\delta ) \|_q + 2 K {\varepsilon} $$
$$ \le 3 K \varepsilon . $$
By Brouwer's fixed point theorem there exists a $\hat \delta $ with
$\| \hat \delta \|_q \le 3 K \varepsilon$ such that
$$ F (\hat \delta ) = \hat \delta .$$
But since $\hat \theta $ is unique we must have
$\hat \theta - \theta^0 = \hat \delta $.
\hfill $\sqcup \mkern -12mu \sqcap$

\subsection{Proof for Section \ref{irrepresentable.section}:
The irrepresentable condition and rates in sup-norm} 

{\bf Proof of Theorem \ref{varselection3.theorem}.} Throughout the proof we assume we are on
${\cal T}_w  \cap {\cal T}_{\rm rem} $.
Recall that for a vector $\theta \in \R^p$, the vector $\theta_S$ is either the
$|S|$-dimensional vector $\{ \theta_j \}_{j \in S}$ or the $p$-dimensional
vector $\{ \theta_j {\rm l} \{ j \in S \} \} $, whichever is appropriate. 
Recall also that $\tilde \theta_S $ is a solution of 
$$ P_n (\dot \rho_{\tilde \theta_S})_S +   \lambda \tilde z_S =0 ,$$
where $\tilde z_{j,S} = {\rm sign} (\tilde \theta_{j,S} )$ if $\tilde \theta_{j,S} \not= 0$, $j \in S\backslash
{\cal J} $,  $\tilde z_{j,S} = 0 $ if $j \in {\cal J}$ and 
$\| \tilde z_S \|_{\infty} \le 1 $. 

We have
$$-w_S + {\cal I}_{S,S} (\tilde \theta_S - \theta_S^0) +
{\cal I}_{S,S^c} \theta_{S^c}^0 -{\rm rem}_S (\tilde \theta_S - \theta^0) + \lambda \tilde z_S =0 . $$
In other words
$$\tilde \theta_S - \theta_S^0= {\cal I}_{S,S}^{-1}  ( {\rm rem} (\tilde \theta_S - \theta^0)) +
v (\tilde \theta_S - \theta_S^0) +u_0 , $$
where
$$v(  \theta_S - \theta_S^0)=  {\cal I}_{S,S}^{-1} ( w_S -\lambda  z_S ), $$
with $z_{j,S} = {\rm sign} (\theta_{j,S})$ if
$\theta_{j,S} \not= 0 $, $ j \in S \backslash {\cal J} $, $z_{j, S} =0 $, $j \in {\cal J} $,  $\| z_S \|_{\infty} \le 1 $, 
and
$$
u_0 = -\Gamma_{S^c}  \theta_{S^c}^0 .$$
Hence
$$\sup_{\theta_S} \| v (\theta_S - \theta_S^0 ) \|_{\infty} \le K_S (\lambda_{w} + \lambda ). $$
and
$$\| u_0 \|_{\infty}  
\le \vvert \Gamma_{S^c } \vvert_1 \lambda_{S^c} \le K_S  \lambda_{S^c}  
 \le K_S \lambda \le K_S (\lambda + \lambda_w )  . $$
 
Let $\delta_S$ be a vector satisfying $\| \delta_S \|_{\infty} \le 3K_S(\lambda_w + \lambda )$.
Since $\lambda_{S^c} \le K_S (\lambda_w + \lambda) $, 
we also have $\| \delta_S - \theta_{S^c}^0 \|_{\infty} \le 
3K_S(\lambda_w + \lambda) $
and so 
$$\| {\rm rem} ( \delta_S + \theta_{S^c}^0 ) \|_{\infty} \le \lambda_{\rm rem} \le \lambda_w + \lambda . $$
By Lemma \ref{Brouwer.lemma} (with $G = {\cal I}_{S,S}^{-1}$ and with
${\rm rem} (\cdot):= {\overline {\rm rem}} (\cdot) $  in
Lemma \ref{Brouwer.lemma} now taken as
${\overline {\rm rem }} (\delta_S)= {\rm rem}( \delta_S  + \theta_{S^c}^0)$) we may now conclude that
$$\| \tilde \theta_S - \theta^0 \|_{\infty} \le 3 K_S (\lambda_w + \lambda) . $$

We now have
$$\tilde \theta_S - \theta_S^0= {\cal I}^{-1} \biggl ( w_S - \lambda \tilde z_S  
 -( {\rm rem} ( \tilde \theta_S - \theta_S^0 ) )_S\biggr ) - \Gamma_{S^c}\theta_{S^c}^0
$$
where
$$\| {\rm rem} ( \tilde \theta_S - \theta^0) \|_{\infty} \le \lambda_{\rm rem} . $$
Thus
$$P_n (\dot \rho_{\tilde \theta_S} )_{S^c} 
= - w_{S^c} + {\cal I}_{S^c, S} (\tilde \theta_S - \theta_S^0) +
{\cal I}_{S,S} \theta_{S^c}^0 -{\rm rem}_{S^c} (\tilde \theta_S - \theta^0) $$
$$ = -C_{S^c}^Tw- \Gamma_{S^c}^T \biggl (   \lambda \tilde z_S  
+ ({\rm rem} ( \tilde \theta_S - \theta_S^0 ))_{S^c}  \biggr )  - {\cal I}_{S^c, S} \Gamma_{S^c}\theta_{S_0^c} +
{\cal I}_{S,S} \theta_{S^c}^0 -({\rm rem} (\tilde \theta_S - \theta^0))_{S^c}  .$$
It follows that
$$\| P_n (\dot \rho_{\tilde \theta_S} )_{S^c} \|_{\infty} \le
\lambda_w (1+ \vvert \Gamma_{S^c} \vvert_1)  + \lambda_{\rm rem} \vvert C_{S^c} \vvert_1  + \vvert \Gamma_{S^c} \vvert_{1, {\rm off}} \lambda +
\kappa_{S^c}  \lambda_{\rm thres} \le \lambda . $$
Therefore, if we define
$$\lambda \tilde z :=  -P_n (\dot \rho_{\tilde \theta_S} )_{S^c}$$
we get that $\tilde \theta_S$ is the solution of
$$P_n \dot \rho_{\tilde \theta_S} + \lambda \tilde z =0 , $$
where $\tilde z_{j} = {\rm sign} (\tilde \theta_{j,S} )$ if $\tilde \theta_{j,S} \not= 0$, $j \notin {\cal J} $, 
$\tilde z_{j,S} =0 $, $j \in {\cal J} $, 
and $\| \hat z \|_{\infty} \le 1 $, i.e.\ $\tilde \theta_S$ is a solution of
(\ref{KKT2.equation}). 
Since the solution is unique we must have $\hat \theta = \tilde \theta_S $.

\hfill $\sqcup \mkern -12mu \sqcap$

\subsection{Proof for Section \ref{precision.section}: Estimating a precision matrix}

{\bf Proof of Lemma \ref{graphicalrem.lemma}.} Let $\Delta := \Theta_S - \Theta_0$. 
It holds that
$$ {\rm rem}  (\Delta)  = ( \Theta_0^{-1} \Delta )^2 (I+ \Theta_0^{-1}\Delta )^{-1} \Theta_0^{-1} . $$
But
$$\vvert (I+ \Theta_0^{-1} \Delta )^{-1} \vvert_{1} \le \sum_{m=0}^{\infty} \vvert ( \Theta_0^{-1} \Delta )^m \vvert_1  $$
and
$$ \vvert ( \Theta_0^{-1} \Delta )^m \vvert_1\le  \vvert \Theta_0^{-1} \Delta  \vvert_1^m
\le (  \vvert \Theta_0^{-1} \vvert_1 \vvert \Delta  \vvert_1 )^m . $$
We have
$$ \vvert \Theta_S - \Theta_0 \vvert_1 \le  \epsilon_0   \vvert  \Theta_S^0 \vvert_0 + \vvert \Theta_{S^c}^0 \vvert_1 . $$
It follows that 
$$ \vvert \Theta_0^{-1} \vvert_1 \vvert \Theta_S - \Theta_0   \vvert_1 \le 
\vvert \Theta_0^{-1} \vvert_1 \biggl ( \varepsilon_0   \vvert  \Theta_S^0 \vvert_0 + \vvert\Theta_{S^c}^0 
\vvert_1 \biggr ) \le \eta_0 $$
and so
$$ \vvert ( I - \Theta_0^{-1} \Delta ) \vvert_1 \le 1/ (1- \eta_0) . $$
We moreover have
$$ \| e_j^T (\Theta_0^{-1} (\Theta_S - \Theta_0 ))^2  \|_1 \le \| e_j^T \Theta_0^{-1} 
(\Theta_S - \Theta_0)  \Theta_0^{-1} \|_1 \varepsilon_0 \le
\vvert \Theta_0^{-1} (\Theta_S - \Theta_0)  \Theta_0^{-1} \vvert_1 \varepsilon_0  $$ $$\le
(\varepsilon_0 \vvert  \Theta_S^0 \vvert_0 + \vvert \Theta_{S^c}^0 \vvert_1)
\vvert \Theta_0^{-1} \vvert_1^2 \epsilon_0 \le\eta_0 \epsilon_0  \vvert \Theta_0^{-1} \vvert_1  .$$
Also
$$
\| (I+ \Theta_0^{-1}(\Theta_S - \Theta_0) )^{-1} \Theta_0^{-1} e_k \|_{\infty} \le
\vvert (I+ \Theta_0^{-1}(\Theta_S - \Theta_0) )^{-1} \vert_1 \vvert \Theta_0^{-1} \vvert_1 $$
$$ \le \vvert \Theta_0^{-1} \vvert_1/(1- \eta_0)$$
So we find
$$\| {\rm rem} (\Delta) \|_{\infty} \le \eta_0 \epsilon_0 \vvert \Theta_0^{-1} \vvert_1^2/(1- \eta_0) . $$
\hfill $\sqcup \mkern -12mu \sqcap$

\bibliographystyle{plainnat}
\bibliography{reference}

\end{document}